\documentclass[a4paper, 11pt]{article}
\usepackage[utf8]{inputenc} 
\usepackage[T1]{fontenc}
\usepackage[english]{babel}
\usepackage{amsmath, amssymb, amsthm}
\usepackage[margin=1in]{geometry}
\usepackage{enumerate}
\usepackage{array}
\usepackage{graphicx}
\usepackage{tikz}
\usepackage{xcolor}
\usepackage{hyperref}
\usepackage{float}
\usepackage{booktabs}
\newtheorem{theorem}{Theorem}[section]
\newtheorem{theoremAB}{Theorem}

\newtheorem{proposition}{Proposition}[section]

\newtheorem{definition}{Definition}[section]
\newtheorem{notation}{Notation}[section]
\newtheorem{lemma}{Lemma}[section]
\newtheorem{observation}{Observation}[section]
\usepackage{listings}
\usepackage{svg}
\usepackage{xcolor}

\lstdefinestyle{pythonstyle}{
    language=Python,
    backgroundcolor=\color{gray!10},
    basicstyle=\ttfamily\small,
    keywordstyle=\color{blue},
    commentstyle=\color{green!50!black},
    stringstyle=\color{orange},
    numbers=left,
    numberstyle=\tiny\color{gray},
    breaklines=true,
    frame=single,
    rulecolor=\color{black},
    tabsize=4,
    showstringspaces=false
}
\lstdefinestyle{consolestyle}{
    backgroundcolor=\color{black!5},
    basicstyle=\ttfamily\small,
    frame=single,
    rulecolor=\color{black},
    breaklines=true,
    showstringspaces=false
}
\usepackage[version=4]{mhchem}
\title{Reverse isoperimetric problem for the center of mass of $\lambda$-convex arcs in $\mathbb{R}^2$}
\author{Sergi Ivars Galiana}
\date{\today}

\begin{document}

\maketitle

\begin{abstract}
This paper focuses on the study of a reverse isoperimetric problem, regarding the minimization of the height of the center of mass of planar curves with constant linear mass density and uniform lower bound of $\lambda$ on their curvatures. We consider several geometric restrictions on the curves and on their curvature, such as imposing $\lambda$-convexity, which constitutes a usual restriction for other reverse isoperimetric problems. We refer to such curves as '$\lambda$-convex standing arcs' which we assume to be placed in the upper half-plane, with endpoints on the horizontal axis. We obtain two results for the minimizer of the height of the center of mass among $\lambda$-convex arcs of a fixed length depending on the restrictions on the position of the endpoints. We show that if the contact angles at the endpoints are bounded above by $\pi/2$, then the minimizer is a circular arc of curvature $\lambda$. Without the angle restriction, the minimizer is a $\lambda$-convex lens, i.e., a closed curve composed of two equal circular arcs of curvature $\lambda$.
\end{abstract}

\section{Introduction and Motivation}

The starting point of this work is the classical isoperimetric problem. Among all simple closed planar curves with fixed length, the circle is the one enclosing the maximal area. More precisely, if $\gamma\subset\mathbb{R}^2$ is a simple closed curve with length $L$ and enclosed area $A$, then
\begin{equation*}
    4\pi A\le L^2,
\end{equation*}
with equality holding if and only if $\gamma$ is a circle.

This problem has been extended to a variety of contexts, such as spaces of constant curvature. A natural question that arises from this classical result is the converse optimization problem: whether there exists a simple closed planar curve that minimizes the area given a fixed perimeter.

Without additional assumptions, the answer is trivial: one can flatten the curve so that it resembles a straight line, making the area arbitrarily small. To yield meaningful and rich geometric results, we can formulate the same minimization question under explicit geometric or curvature constraints. A widely discussed approach has been to impose a lower bound on the curvature of the curve, that is, requiring its curvature to be at least $\lambda$ at every point, for some fixed $\lambda>0$. This condition is formalized through the concept of $\lambda$-convexity, which will appear repeatedly throughout this work and will be formally defined later in Definition \ref{def lambda-convex}.

In this direction, a sharp reverse isoperimetric inequality under a lower curvature bound was first introduced by Keith Ball in \cite{ball} and \cite{ball2}. Ball proved that among all simple closed curves $\gamma\subset\mathbb{R}^2$ with fixed length $L$ whose curvature is everywhere at least $\lambda$, the minimizer is the boundary of the symmetric intersection of two disks of radius $1/\lambda$. This extremal body is known as the $\lambda$-convex lens (see Figure \ref{fig:convex lens}). More precisely, the area enclosed by any such curve is bounded from below by the area of the $\lambda$-convex lens of equal length. The uniqueness of the minimizer was later rigorously established by Franck Barthe in \cite{barthe}, and the complete statement was proven in \cite{kostya3}, with an alternative geometric approach given in \cite{RIPn=2}.

Extending the reverse isoperimetric problem for $\lambda$-convex curves to higher dimensions becomes increasingly difficult. The generalization of this problem was finally solved in dimension three by Drach and Tatarko \cite{kostya1}, who proved that the $\lambda$-convex lens remains the unique volume minimizer among all $\lambda$-convex bodies of a given surface area. For dimensions greater than or equal to four, the problem was generalized and solved just recently by Drach and Tatarko \cite{drach2026completesolutionreverseisoperimetric}.

Reverse isoperimetric problems can also go beyond the classical interplay between volume and surface area. Similar tools have been applied to optimize other fundamental geometric measurements under curvature constraints, such as the \textit{inradius} of a convex domain, defined as the radius of the largest ball contained in it. This problem was first solved in the two-dimensional case by Milka \cite{inradius1}, and later extended to higher dimensions, first for bodies with a prescribed volume by Bezdek \cite{inradius2}, and then for bodies with a prescribed surface area in \cite{kostya1}. For other examples of related reverse isoperimetric problems, see \cite{example1}.

On the other hand, a parallel branch of the literature explores reverse isoperimetric problems under upper curvature bounds, formalized through the concept of $\lambda$-concavity. In this setting, the minimization of volume among $\lambda$-concave isoperimetric domains yields an entirely different extremal shape: the $\lambda$-sausage, i.e., the convex hull of two balls of radius $1/\lambda$. This problem was fully resolved in \cite{kostya2concave} (see \cite{alternativeconcave} for an alternative proof).

While the aforementioned literature heavily centers on optimizing global geometric measurements such as area, volume, surface area, or inradius under curvature constraints, the focus of this work shifts towards a fundamentally different geometric quantity: the spatial distribution of mass along a curve. Specifically, we investigate the height of the center of mass of a geometric curve endowed with a uniform linear density.

To formalize the setting of our study, we restrict our attention to arcs of rectifiable curves embedded in the upper half-plane whose endpoints lie on the $x$-axis; throughout this work, we will refer to such arcs as \textbf{standing arcs}. We propose and analyze the following optimization problem:

\begin{center}
    \textit{Among all standing arcs of a fixed length $L$ whose curvature is everywhere bounded below by a constant $\lambda>0$, which geometrical configuration minimizes the vertical coordinate of its center of mass?}
\end{center}

Much like the classical area-minimization problems, the lower curvature bound is strictly necessary to prevent the trivial case where the standing arc is indefinitely flattened against the $x$-axis, bringing its center of mass arbitrarily close to zero. Intuitively, one might hypothesize that a smooth circular arc of radius $1/\lambda$ would emerge as the absolute minimizer for this functional. However, this intuition fails in general; a counterexample is provided by the $\lambda$-convex lens.

The primary objective of this work is to analyze this optimization problem in detail and to characterize the shapes that minimize the height of the center of mass under different geometric assumptions.

\section{Basic definitions and statement of the main results}

We begin by introducing the notation and basic concepts necessary for the statement of the main results and that will also be used throughout the article.

\begin{definition}
    Let $\gamma\subset\mathbb{R}^2$ be a rectifiable curve with linear continuous mass density $\rho(x,y)$ that is not identically zero. The coordinates of the \textbf{center of mass} of $\gamma$ are defined by
    \begin{gather*}
        x_{CM}=\frac{\int_{\gamma}x\rho(x,y)\,\text{d}L}{\int_{\gamma}\rho(x,y)\,\text{d}L}, \\
        y_{CM}=\frac{\int_{\gamma}y\rho(x,y)\,\text{d}L}{\int_{\gamma}\rho(x,y)\,\text{d}L},
    \end{gather*}
    where $\int_\gamma dL$ represents integration with respect to arc length along $\gamma$. In particular, the denominators refer to the total mass of the curve.
\end{definition}

    When referring to the length and the height of the center of mass ($y_{CM}$) of a given curve $\gamma$ we will write $L(\gamma$) and $CM(\gamma$), respectively.

\begin{observation}
    Throughout this work we restrict to curves with constant linear density. In this case our expression simplifies to:
    \begin{equation}
    \label{center of mass eqn}
        (x_{CM}, y_{CM})=\left(\frac{\int_{\gamma}x\,\text{d}L}{\int_{\gamma}\text{d}L}, \frac{\int_{\gamma}y\,\text{d}L}{\int_{\gamma}\text{d}L}\right)=\left(\frac{\int_{\gamma}x\,\text{d}L}{L(\gamma)}, \frac{\int_{\gamma}y\,\text{d}L}{L(\gamma)}\right).
    \end{equation}
\end{observation}

All the results of the paper take place in the upper half plane of $\mathbb{R}^2$, which can be explicitly defined as: $(\mathbb{R}^2)^+:=\{(x,y)\in\mathbb{R}^2 \ | \ y\ge0\}$. We can now define the following:

\begin{definition}
Given a standing arc $\gamma\subset(\mathbb{R}^2)^+$, we define its set of \textbf{endpoints} as $E=\{(x,y)\in\mathbb{R}^2|y=0\}\cap\gamma$. We denote by $|E|$ the number of points of $E$.
\end{definition}
\begin{observation}
    We will only consider curves that have at least one endpoint ($|E|\ge1$). This makes sense in the context of our study, since if a curve does not have endpoints as we just described them, then by performing a vertical translation downward until it has at least one (while still being contained in $(\mathbb{R}^2)^+$) we obtain a new curve with a strictly lower center of mass but with the same length and curvature.
\end{observation}
For two-dimensional spaces, a basic definition of $\lambda$-convexity is the following:

\begin{definition}
\label{def 1.2}
    Let $K\subset\mathbb{R}^2$ be a convex body with boundary $\partial K=\gamma$. We say that $K$ is a $\lambda$-convex body (and $\gamma$ is a $\lambda$-convex curve) if there exists $\lambda>0$ such that the curvature $k$ of $\gamma$ satisfies:
    \begin{equation*}
        k(p)\ge\lambda \ \ \ \ \forall p\in\gamma.
    \end{equation*}
\end{definition}

When we say that a standing arc $\gamma$ is $\lambda$-convex, we mean that $\gamma$ together with the segment of the $x$-axis joining its endpoints bounds a convex domain $D$, and that $\gamma$ has curvature greater than or equal to $\lambda$ at every point.

\begin{observation}
    Given these restrictions, we have that $|E|\in\{1,2\}$. Furthermore, $|E|=1$ if and only if the curve begins and ends at the same point in the $x$-axis.
\end{observation}

For the following definition, we need to work with the tangent of $\gamma$ at its endpoints. In upcoming sections we will see that, under certain approximations, it is possible to assume that the tangent of a standing arc is defined everywhere; however, we do not require this assumption just yet. For our immediate purposes, the tangent at the right-most endpoint refers to the left-hand tangent (and vice versa for the left-most endpoint), both of which are always well-defined.

\begin{definition}
\label{contact angles}
    Let $\gamma$ be a standing arc and let $p\in E$ be an endpoint. The \textbf{contact angle} at $p$ is the angle between the tangent line to $\gamma$ at $p$ and the $x$-axis.

    If $\gamma$ has two endpoints, we measure this angle clockwise from the tangent at the left-most endpoint and counterclockwise at the right-most endpoint. These angles are denoted by $\theta_1(\gamma)$ and $\theta_2(\gamma)$, respectively.

    If $\gamma$ has only one endpoint, $\theta_1(\gamma)$ and $\theta_2(\gamma)$ are measured as before, but using the left and right supporting tangents at $p$, respectively.
\end{definition}
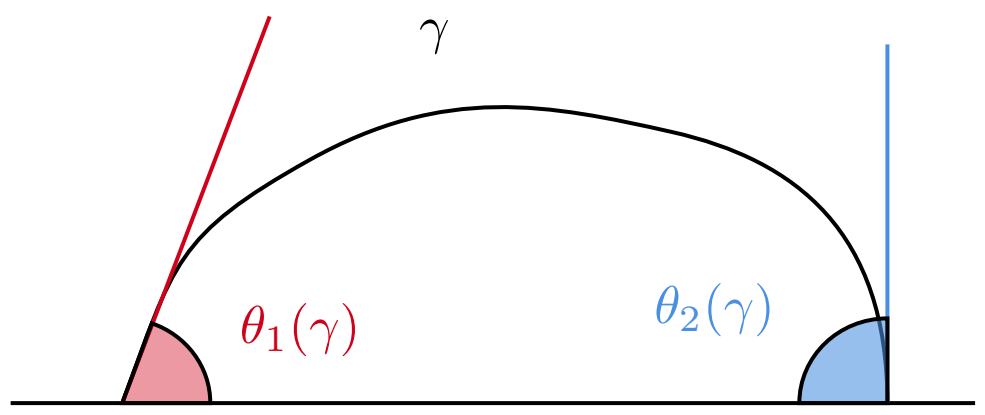
\begin{figure}[h!]
    \centering

\tikzset{every picture/.style={line width=0.75pt}} 

\begin{tikzpicture}[x=0.75pt,y=0.75pt,yscale=-1,xscale=1]

\draw    (99,130.6) -- (341,130.6) ;
\draw    (127,130.6) .. controls (140,95.6) and (140,88.6) .. (174,69.6) .. controls (208,50.6) and (234,55.6) .. (265,62.6) .. controls (296,69.6) and (319,87.6) .. (319,130.6) ;
\draw [color={rgb, 255:red, 208; green, 2; blue, 27 }  ,draw opacity=1 ]   (127,130.6) -- (164,33.6) ;
\draw  [fill={rgb, 255:red, 208; green, 2; blue, 27 }  ,fill opacity=0.4 ] (134.42,110.58) .. controls (142.97,113.51) and (149.1,121.37) .. (149.1,130.6) -- (127,130.6) -- cycle ;
\draw [color={rgb, 255:red, 74; green, 144; blue, 226 }  ,draw opacity=1 ]   (319,40.6) -- (319,130.6) ;
\draw  [fill={rgb, 255:red, 74; green, 144; blue, 226 }  ,fill opacity=0.58 ] (296.9,130.6) .. controls (296.9,118.86) and (306.79,109.35) .. (319,109.35) -- (319,130.6) -- cycle ;

\draw (200,33) node [anchor=north west][inner sep=0.75pt]  [font=\small] [align=left] {$\displaystyle \gamma $};
\draw (155,104) node [anchor=north west][inner sep=0.75pt]  [font=\small,color={rgb, 255:red, 208; green, 2; blue, 27 }  ,opacity=1 ] [align=left] {$\displaystyle \theta _{1}( \gamma )$};
\draw (259,99) node [anchor=north west][inner sep=0.75pt]  [font=\small,color={rgb, 255:red, 74; green, 144; blue, 226 }  ,opacity=1 ] [align=left] {$\displaystyle \theta _{2}( \gamma )$};

\end{tikzpicture}
\caption{Representation of the concept of contact angle.}
\end{figure}

A particular type of $\lambda$-convex curve of special interest, not only in our study but also in other reverse isoperimetric inequalities with curvature restrictions, is the $\lambda$-convex lens.

\begin{figure}[h]
    \centering

\tikzset{every picture/.style={line width=0.75pt}} 

\begin{tikzpicture}[x=0.75pt,y=0.75pt,yscale=-1,xscale=1]

\draw    (99,138.6) -- (384,138.6) ;
\draw  [draw opacity=0] (371.79,98.5) .. controls (344.43,122.47) and (298.25,138.23) .. (245.9,138.23) .. controls (193.09,138.23) and (146.57,122.2) .. (119.3,97.88) -- (245.9,47.66) -- cycle ; \draw   (371.79,98.5) .. controls (344.43,122.47) and (298.25,138.23) .. (245.9,138.23) .. controls (193.09,138.23) and (146.57,122.2) .. (119.3,97.88) ;  
\draw  [draw opacity=0] (119.3,97.88) .. controls (146.63,73.87) and (192.79,58.05) .. (245.14,57.97) .. controls (297.94,57.88) and (344.48,73.83) .. (371.79,98.1) -- (245.27,148.52) -- cycle ; \draw   (119.3,97.88) .. controls (146.63,73.87) and (192.79,58.05) .. (245.14,57.97) .. controls (297.94,57.88) and (344.48,73.83) .. (371.79,98.1) ;  
\draw [color={rgb, 255:red, 74; green, 144; blue, 226 }  ,draw opacity=1 ] [dash pattern={on 4.5pt off 4.5pt}]  (119.3,97.88) -- (371.79,98.5) ;

\end{tikzpicture}
\caption{Example of a $\lambda$-convex lens and its main axis (dotted blue line).}
\label{fig:convex lens}
\end{figure}
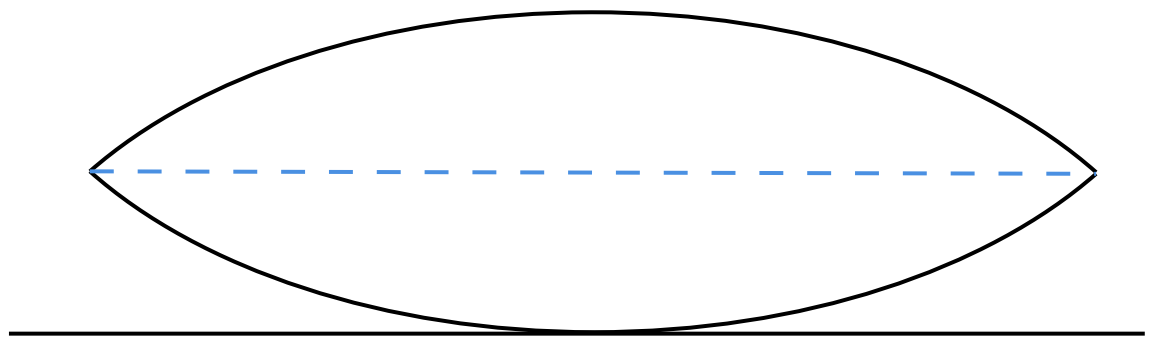

\begin{definition}
\label{def lambda-convex}
    A \textbf{$\lambda$-convex lens} is a closed curve composed of two circular arcs of radius $1/\lambda$ having the same length and meeting at their endpoints. The segment joining its 2 endpoints will be referred to as its \textbf{major axis}. See Figure \ref{fig:convex lens}.
\end{definition}

\subsection{Reverse isoperimetric inequalities for $\lambda$-convex standing arcs}

The two main results of this work are the following:

\begin{theoremAB}
\label{theorem A}
    Let $\gamma$ be a standing $\lambda$-convex arc in $(\mathbb{R}^2)^+$ with length $L$. Then
    \begin{equation*}
        CM(\gamma)\ge CM(\gamma^*)=\frac{1}{\lambda}\left(1-\cos\left(\frac{L\lambda}{4}\right)\right),
    \end{equation*}
    where $\gamma^*$ is the standing arc with the geometry of the $\lambda$-convex lens of length $L$ whose major axis is parallel to the x-axis.
\end{theoremAB}

Furthermore, if we add an extra condition regarding the size of the contact angles we obtain a different result:

\begin{theoremAB}
\label{theorem B}
    Let $\gamma$ be a standing arc satisfying the assumptions of Theorem \ref{theorem A}. Assume in addition that the contact angles are at most $\pi/2$ $(\theta_i(\gamma)\le\pi/2$ \ $\forall i)$. Then,
    \begin{equation*}
        CM(\gamma)\ge CM(\tilde\gamma)=\frac{1}{\lambda}\left(\sin\left(\frac{L\lambda}{2}\right)\frac{2}{L\lambda}-\cos\left(\frac{L\lambda}{2}\right)\right),
    \end{equation*}
    where $\tilde\gamma$ is the standing circular arc of radius 1/$\lambda$ and length $L$.
\end{theoremAB}

The strategy for proving these results is as follows: we first establish Theorem \ref{theorem B} by exploiting the arc-length parametrization and its intrinsic relationship with the curvature. For Theorem \ref{theorem A} we first reduce the problem to vertically symmetric cases, and then apply the hinge technique (see \cite{burago2001metric} and  \cite{petersen2006riemannian}) to lower the curve's center of mass until obtaining a $\lambda$-convex lens.

To prove these statements, we must initially assume certain regularity properties of the curve, such as requiring it to be differentiable everywhere or to possess piecewise constant curvature. To achieve this, we will employ approximation methods for $\lambda$-convex curves such as the curved polygon approximation; which has previously been studied in the literature of the reverse isoperimetric problem \cite{kostya1}.

\section{Preliminaries}
\subsection{Foundations and Notation of $\lambda$-convex standing arcs}
In this section we collect some definitions and results that will be used later. 

We have observed that analyzing symmetric configurations significantly simplifies the proof, particularly for Theorem \ref{theorem A}. For this reason, the purpose of the following definitions is to construct, from an arbitrary standing arc $\gamma$, a symmetric standing arc of equal length whose center of mass is lower than or at the same height as that of the original curve.

\begin{figure}[h!]
\label{hhp lines image}
    \centering

\tikzset{every picture/.style={line width=0.75pt}} 

\begin{tikzpicture}[x=0.55pt,y=0.55pt,yscale=-1,xscale=1]

\draw    (70.45,179.92) .. controls (95,128.6) and (93,130.6) .. (108,115.6) .. controls (123,100.6) and (131.22,93.02) .. (149.61,84.31) .. controls (168,75.6) and (178,70.6) .. (212,89.6) .. controls (246,108.6) and (238,102.6) .. (259.85,122.47) .. controls (281.71,142.34) and (290.06,164.32) .. (293.2,180.32) ;
\draw    (47.2,179.6) -- (323.87,180.27) ;
\draw [color={rgb, 255:red, 208; green, 2; blue, 27 }  ,draw opacity=1 ] [dash pattern={on 4.5pt off 4.5pt}]  (63.53,126.85) -- (316.49,126.85) ;
\draw    (323.87,180.27) -- (647.2,179.6) ;
\draw [color={rgb, 255:red, 74; green, 144; blue, 226 }  ,draw opacity=1 ] [dash pattern={on 4.5pt off 4.5pt}]  (496.96,59.87) -- (496.56,202.67) ;
\draw    (390.45,179.92) .. controls (415,128.6) and (413,130.6) .. (428,115.6) .. controls (443,100.6) and (451.22,93.02) .. (469.61,84.31) .. controls (488,75.6) and (498,70.6) .. (532,89.6) .. controls (566,108.6) and (558,102.6) .. (579.85,122.47) .. controls (601.71,142.34) and (610.06,164.32) .. (613.2,180.32) ;

\draw (56.75,143.8) node [anchor=north west][inner sep=0.75pt]  [font=\small] [align=left] {$\displaystyle \gamma _{-}$};
\draw (264.93,104.18) node [anchor=north west][inner sep=0.75pt]  [font=\footnotesize,color={rgb, 255:red, 208; green, 2; blue, 27 }  ,opacity=1 ] [align=left] {$\displaystyle HHP( \gamma )$};
\draw (503.2,56.6) node [anchor=north west][inner sep=0.75pt]  [font=\footnotesize,color={rgb, 255:red, 74; green, 144; blue, 226 }  ,opacity=1 ] [align=left] {$\displaystyle HVP( \gamma )$};
\draw (296.75,143.8) node [anchor=north west][inner sep=0.75pt]  [font=\small] [align=left] {$\displaystyle \gamma _{-}$};
\draw (114.75,67.8) node [anchor=north west][inner sep=0.75pt]  [font=\small] [align=left] {$\displaystyle \gamma _{+}$};
\draw (408.75,90.8) node [anchor=north west][inner sep=0.75pt]  [font=\small] [align=left] {$\displaystyle \gamma _{-}$};
\draw (580.75,90.8) node [anchor=north west][inner sep=0.75pt]  [font=\small] [align=left] {$\displaystyle \gamma _{+}$};

\end{tikzpicture}
    \caption{Illustration of the concepts $HHP(\gamma)$ and $HVP(\gamma)$, respectively.}
    \label{fig:enter-label}
\end{figure}
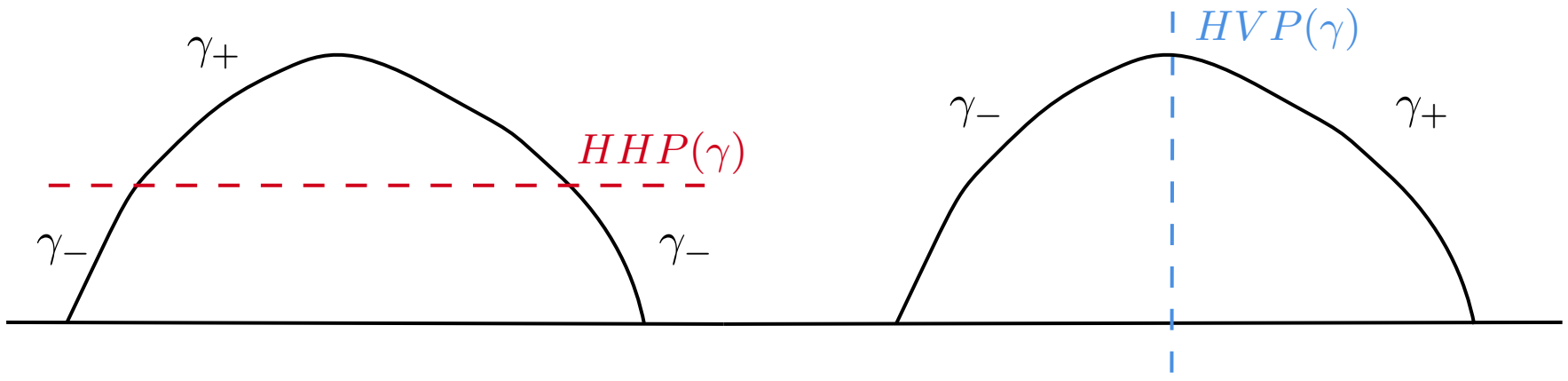

\begin{definition}
\label{def half horizontal}
    Let $\gamma\subset(\mathbb{R}^2)^+$ be a $\lambda$-convex standing arc. The \textbf{half-perimeter horizontal line}, \textbf{HHP($\gamma$)}, is the line parallel to the x-axis that divides $\gamma$ into two parts of equal length. See Figure \ref{hhp lines image}.
\end{definition}

\begin{definition}
\label{def half vertical}
    Let $\gamma\subset(\mathbb{R}^2)^+$ be a $\lambda$-convex standing arc. The \textbf{half-perimeter vertical line}, \textbf{HVP($\gamma$)}, is the line parallel to the y-axis that divides $\gamma$ into two parts of equal length. See Figure \ref{hhp lines image}.
\end{definition}

We can justify the existence of $HHP(\gamma)$ and $HVP(\gamma)$ as follows: fix a strictly convex standing arc $\gamma$ and let $r_H$ be a horizontal line. Move $r_H$ continuously upwards (from a position strictly below $\gamma$). Then, the length of the part of $\gamma$ lying below $r_H$ varies continuously from 0 to $L(\gamma)$; since the function $f(h)$ that returns the arc length of $\gamma$ below $y=h$ is continuous and strictly increasing. Hence, by the Intermediate Value Theorem there exists a position of $r_H$ for which this length equals $L(\gamma)/2$, i.e., $r_H=HHP(\gamma)$, proving thus the existence of the half horizontal perimeter line. The vertical case is analogous.

\begin{notation}
    When dividing a standing arc $\gamma$ with the $HHP(\gamma)$ we will denote by $\gamma_-$ and $\gamma_+$ the portions of the curve lying below and above this line, respectively. Since we will not draw $HHP(\gamma)$ and $HVP(\gamma)$ simultaneously at any point in the paper we also denoted by $\gamma_-$ and $\gamma_+$ the portions of the curve lying to the left and to the right of $HVP(\gamma)$, respectively. See Figure \ref{hhp lines image}.
\end{notation}

In our context, it will be convenient to use the following characterization of $\lambda$-convexity, which is equivalent to Definition \ref{def 1.2}.
\begin{lemma}
\label{def convex 2}
    Let $\gamma\subset\mathbb{R}^2$ be a simple closed curve. Then $\gamma$ is $\textbf{$\lambda$-convex}$ if for every p $\in\gamma$ there exists a neighborhood, $\mathcal{U}_p\subset\mathbb{R}^2$, of p and a closed disk of radius $1/\lambda$, $B_p(1/\lambda)\subset\mathbb{R}^2$,  such that: $p\in\partial B_p(1/\lambda)$ and $\mathcal{U}_p\cap\gamma\subset B_p(1/\lambda)$.
\end{lemma}

See \cite{kostya1} for a generalized definition in any dimension and for a visual example. This lemma is a direct consequence of a classical result for $\lambda$-convex curves that allows us to bypass standard differentiability assumptions (a point we discuss further in the following subsections): Blaschke's Rolling Theorem.

\begin{theorem}[Blaschke's Rolling Theorem]
\label{blaschke}
Let D $\subset\mathbb{R}^2$ be a convex domain with a $\lambda$-convex boundary $\gamma=\partial D$. Then, for every $p\in\gamma$, there exists a closed disk $B_p(1/\lambda)$ of radius 1/$\lambda$ is such that $p\in\partial B_p(1/\lambda)$ then $D\subset B_p(1/\lambda)$.
\end{theorem}

We will refer to the disks $B_p(1/\lambda)$ as \textbf{supporting disks}.

This well-known result in differential geometry has been stated and proven in multiple contexts. See, for example, \cite{blaschke2dim} for the first proof of the theorem for domains with smooth boundaries in $\mathbb{R}^2$ and \cite{blaschke} or \cite{blaschke2} for the general statement in $\mathbb{R}^n$. Furthermore, see \cite{blaschke3} for a new generalization of the theorem for convex domains in Riemannian manifolds of bounded sectional curvature and arbitrary dimension.

The following definition is formulated under the assumption that our standing arcs possess a well-defined tangent at all points. As previously noted, this regularity condition will be rigorously established in upcoming sections by combining Blaschke's Rolling Theorem with other approximation techniques.
\begin{definition}
\label{def section angles}
    Let $\gamma$ be a $\mathcal{C}^1$ $\lambda$-convex standing arc, and let $HVP(\gamma)$ be as in Definition \ref{def half vertical}. Set $p=HVP(\gamma)\cap\gamma$ (choosing the point of greatest height if the intersection is not unique). The \textbf{section angles}, denoted as $\alpha_-$ and $\alpha_+$, are defined as the angles between $HVP(\gamma)$ and the tangent of $\gamma$ at $p$, measured clockwise and counterclockwise from the vertical line, respectively.
\end{definition}

\begin{figure}[h]
    \centering

\tikzset{every picture/.style={line width=0.75pt}} 

\begin{tikzpicture}[x=0.75pt,y=0.75pt,yscale=-1,xscale=1]

\draw    (392.45,179.25) .. controls (431,130) and (424.78,140.83) .. (438.89,122.42) .. controls (453,104) and (485.8,95.03) .. (497.03,94.2) .. controls (508.26,93.38) and (512.52,93) .. (540,109.2) .. controls (567.48,125.4) and (552.15,115.4) .. (577,134.2) .. controls (601.85,153) and (612.06,163.65) .. (615.2,179.65) ;
\draw [color={rgb, 255:red, 74; green, 144; blue, 226 }  ,draw opacity=1 ] [dash pattern={on 4.5pt off 4.5pt}]  (488.96,78.87) -- (488.99,163.2) -- (489,194) ;
\draw    (367,178.6) -- (638,179.6) ;
\draw [color={rgb, 255:red, 208; green, 2; blue, 27 }  ,draw opacity=1 ]   (403,114) -- (554,81.6) ;
\draw  [fill={rgb, 255:red, 208; green, 2; blue, 27 }  ,fill opacity=0.46 ] (489.03,113.2) .. controls (480.9,113.49) and (473.38,107.81) .. (471.41,99.34) -- (488.21,95.45) -- cycle ;
\draw  [fill={rgb, 255:red, 208; green, 2; blue, 27 }  ,fill opacity=0.43 ] (507.14,91.17) .. controls (507.65,93.11) and (507.83,95.15) .. (507.61,97.23) .. controls (506.64,106.28) and (498.59,113) .. (489.03,113.2) -- (488.23,95.17) -- cycle ;
\draw  [color={rgb, 255:red, 0; green, 0; blue, 0 }  ,draw opacity=1 ][fill={rgb, 255:red, 80; green, 227; blue, 194 }  ,fill opacity=1 ] (484.33,95.17) .. controls (484.33,93.02) and (486.07,91.27) .. (488.23,91.27) .. controls (490.38,91.27) and (492.13,93.02) .. (492.13,95.17) .. controls (492.13,97.33) and (490.38,99.07) .. (488.23,99.07) .. controls (486.07,99.07) and (484.33,97.33) .. (484.33,95.17) -- cycle ;

\draw (505.2,153.6) node [anchor=north west][inner sep=0.75pt]  [font=\footnotesize,color={rgb, 255:red, 74; green, 144; blue, 226 }  ,opacity=1 ] [align=left] {$\displaystyle HVP( \gamma )$};
\draw (458,107) node [anchor=north west][inner sep=0.75pt]  [color={rgb, 255:red, 208; green, 2; blue, 27 }  ,opacity=1 ] [align=left] {$\displaystyle \alpha _{-}$};
\draw (513,105) node [anchor=north west][inner sep=0.75pt]  [color={rgb, 255:red, 208; green, 2; blue, 27 }  ,opacity=1 ] [align=left] {$\displaystyle \alpha _{+}$};
\draw (470,72) node [anchor=north west][inner sep=0.75pt]  [font=\normalsize,color={rgb, 255:red, 80; green, 227; blue, 194 }  ,opacity=1 ] [align=left] {$\displaystyle p$};
\draw (383,125) node [anchor=north west][inner sep=0.75pt]   [align=left] {$\displaystyle \gamma _{-}$};
\draw (588,101) node [anchor=north west][inner sep=0.75pt]   [align=left] {$\displaystyle \gamma _{+}$};

\end{tikzpicture}
    \caption{Illustration of the concept of section angles.}
    \label{fig:section angles}
\end{figure}
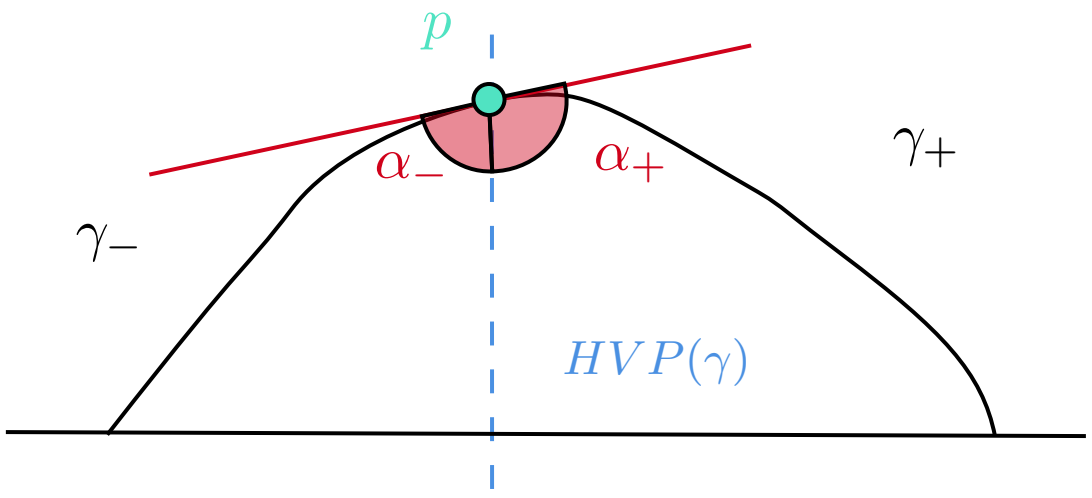

Lastly, we will also need the following result for the proof of Theorem \ref{theorem A}.

\begin{theorem}[The Bent Bow Theorem]
\label{Bent Bow Problem}
    Let $\gamma_1$ and $\gamma_2$ be two standing arcs of the same length $L$. Denote by $k_1(s)$ and $k_2(s)$ their curvatures at the point $s\in[0,L]$ (where $s$ represents the arc length parameter). Assume that $k_1(s)\ge k_2(s)$ $\forall s\in[0,L]$ and $\int_0^Lk_1(s)\text{d}s<\pi$.
    
    Then $\overline{\gamma_1(0)\gamma_1(L)}\le\overline{\gamma_2(0)\gamma_2(L)}$,
    where $\overline{\gamma_i(0)\gamma_i(L)}$, $i=1,2$, denotes the horizontal distance between the endpoints of $\gamma_i$.
\end{theorem}

This is a classical result, whose proof can be found in most introductory textbooks on differential geometry (see, for example, \cite{diffgeombook}).

\subsection{Approximation methods}

\subsubsection{Curved polygon approximation}
\label{subsubsection: curved polygons}

Throughout this section, $\gamma$ will denote a $\lambda$-convex standing arc. To begin, we apply Theorem \ref{blaschke} to a finite set of points on our curve $\gamma$. We now define a new standing arc with the boundary of the region given by the intersection of the supporting disks. 

\begin{definition}
\label{def curved polygons}
    Let $P\subset\gamma$ be a set of $n<\infty$ points. For each $p\in P$, consider the supporting disk $B_p(1/\lambda)$. We define
    \begin{equation*}
    \gamma_n:=\partial\left(\bigcap_{p\in P}B_p(1/\lambda)\right)\cap(\mathbb{R}^2)^+.
    \end{equation*}
    The curves obtained through this procedure will be called \textbf{curved polygons}.  
\end{definition}

The following observation will be useful later:

\begin{observation}
    Let $\{\gamma_i\}$ be a sequence of curved polygons associated with finite sets $P_i$ such that $P_{n}\subset P_m$ whenever $n<m$. Since each curved polygon is an intersection of supporting disks, the domain bounded by $\gamma_m$ is contained in the domain bounded by $\gamma_{n}$. In particular, $L(\gamma_{n})>L(\gamma_m)$ (as can be seen by a direct application of the Bent Bow Theorem \ref{Bent Bow Problem}).
\end{observation}

We now examine the convergence of the sequence of curved polygons $\{\gamma_n\}$ to the curve $\gamma$ as the cardinality of their underlying finite sets grows. To quantify the distance between these curves, we employ the Hausdorff metric.

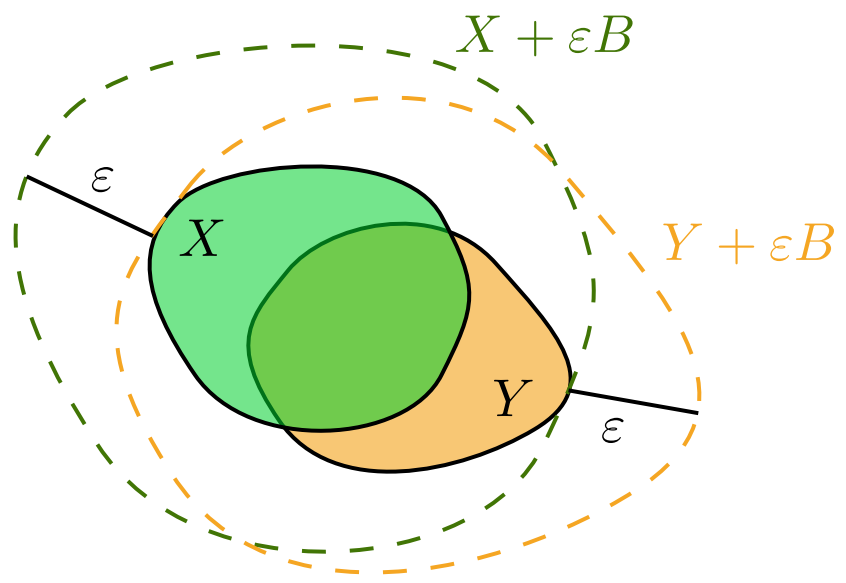
\begin{figure}[h!]
    \centering
\tikzset{every picture/.style={line width=0.75pt}} 

\begin{tikzpicture}[x=0.55pt,y=0.55pt,yscale=-1.2,xscale=1.2]

\draw  [fill={rgb, 255:red, 245; green, 166; blue, 35 }  ,fill opacity=0.63 ] (307.15,134.84) .. controls (317.45,122.12) and (349.16,113.59) .. (366.34,132.84) .. controls (383.52,152.09) and (398.36,167.56) .. (377.43,180.13) .. controls (356.5,192.7) and (320.27,200.4) .. (304.65,177.75) .. controls (289.03,155.11) and (296.84,147.56) .. (307.15,134.84) -- cycle ;
\draw  [fill={rgb, 255:red, 2; green, 208; blue, 46 }  ,fill opacity=0.55 ] (275.75,115.48) .. controls (287.86,103.4) and (339.79,98.87) .. (350.88,119.37) .. controls (361.97,139.86) and (360.87,144.92) .. (350.88,164.66) .. controls (340.88,184.4) and (296.21,187.3) .. (280.6,164.66) .. controls (264.98,142.01) and (263.65,127.56) .. (275.75,115.48) -- cycle ;
\draw  [color={rgb, 255:red, 65; green, 117; blue, 5 }  ,draw opacity=1 ][dash pattern={on 4.5pt off 4.5pt}] (244.08,90.41) .. controls (265.98,67.28) and (359.92,58.61) .. (379.98,97.85) .. controls (400.04,137.09) and (398.06,146.77) .. (379.98,184.57) .. controls (361.89,222.36) and (281.09,227.93) .. (252.84,184.57) .. controls (224.59,141.21) and (222.19,113.53) .. (244.08,90.41) -- cycle ;
\draw  [color={rgb, 255:red, 245; green, 166; blue, 35 }  ,draw opacity=1 ][dash pattern={on 4.5pt off 4.5pt}] (278.7,111.09) .. controls (297.34,86.73) and (354.7,70.4) .. (385.77,107.26) .. controls (416.85,144.11) and (443.69,173.74) .. (405.83,197.8) .. controls (367.97,221.87) and (302.43,236.61) .. (274.18,193.25) .. controls (245.92,149.89) and (260.05,135.44) .. (278.7,111.09) -- cycle ;
\draw    (387.1,169.05) -- (408.95,172.89) -- (424.1,175.55) ;
\draw    (232.6,108.05) -- (268.6,125.05) ;

\draw (363.62,164.33) node [anchor=north west][inner sep=0.75pt]  [font=\small] [align=left] {$\displaystyle Y$};
\draw (274.1,118.66) node [anchor=north west][inner sep=0.75pt]  [font=\small] [align=left] {$\displaystyle X$};
\draw (411.9,119.8) node [anchor=north west][inner sep=0.75pt]  [font=\small,color={rgb, 255:red, 245; green, 166; blue, 35 }  ,opacity=1 ] [align=left] {$\displaystyle Y+\varepsilon B$};
\draw (352.9,60.5) node [anchor=north west][inner sep=0.75pt]  [font=\small,color={rgb, 255:red, 65; green, 117; blue, 5 }  ,opacity=1 ] [align=left] {$\displaystyle X+\varepsilon B$};
\draw (394.6,175.5) node [anchor=north west][inner sep=0.75pt]  [font=\small] [align=left] {$\displaystyle \varepsilon $};
\draw (249.1,104) node [anchor=north west][inner sep=0.75pt]  [font=\small] [align=left] {$\displaystyle \varepsilon $};

\end{tikzpicture}

    \caption{Illustration of the Hausdorff distance between two subsets $X$ and $Y$ in $\mathbb{R}^2$.}
    \label{image hausdorff distance}
\end{figure}

\begin{definition}
    Let $(\mathbb{R}^2,d)$ be a Euclidean space equipped with the standard Euclidean metric. For a pair of compact subsets $X,Y$ $\subseteq$ $\mathbb{R}^2$ we define the \textbf{Hausdorff metric} (see Figure \ref{image hausdorff distance}) as:
    \begin{equation*}
        d_H(X,Y)=\max\{ \sup_{x\in X}d(x,Y),\sup_{y\in Y}d(y,X)\}.
    \end{equation*}
    Another equivalent way to define the Hausdorff metric uses the unit closed disk B. This alternative definition reads as follows:
    \begin{equation*}
        d_H(X,Y)=\inf\{\varepsilon>0 | Y\subseteq X+\varepsilon B \; \text{and} \; X\subseteq Y+\varepsilon B\}.
    \end{equation*}
\end{definition}

Two compact sets $X,Y$ are equal if and only if their Hausdorff distance is 0. Therefore, when we consider sequences of curves converging in the Hausdorff metric we are actually studying the convergence of compact subsets in $\mathbb{R}^2$. With this notion of convergence, we can state the following lemma.

\begin{lemma}
\label{convergence lemma}
    For a sequence of sets of points $\{P_n\}\subset\gamma$ such that $P_n\subset P_{n+1}$ and $\max_{x\in\gamma}d(x,P_n)\rightarrow0$ as $n\rightarrow\infty$, define the curved polygons $\gamma_n$ associated to each of these finite sets of points. Then, the sequence of curved polygons $\{\gamma_n\}_{n\in\mathbb{N}}$ converges to the original $\lambda$-convex curve $\gamma$ in the Hausdorff metric, i.e., $\gamma_n\xrightarrow{H}\gamma$ as $n\rightarrow\infty$.
\end{lemma}

This type of approximation has been used in reverse isoperimetric problems and other areas of convex geometry. For examples of applications of this approximation, see \cite{kostya1} and \cite{polygons}.

\subsubsection{Small circles approximation}
\label{subsubsection: small circles approxiamtion}
When intersecting supporting disks, we might be left with certain points where the derivative is not defined. However, by the construction of our curved polygons, these non-differentiable points are isolated. We can eliminate this lack of differentiability by inserting a small circular arc contained within the domain bounded by the curved polygon. Specifically, if the supporting disks have radius $1/\lambda$, we select an smoothing circle of radius $\varepsilon \ll 1/\lambda$ that is tangent to the boundary of the curved polygon at two points, one on each side of the corner (the non-differentiable point). We also choose $\varepsilon$ small enough to ensure the modified curve still encloses the original curve $\gamma$, as the curved polygon did. This configuration is always realizable due to the continuity of the curve and the positive distance between the two contact points. This construction is illustrated in the figure below.
\begin{figure}[h!]
\centering    

\tikzset{every picture/.style={line width=0.75pt}} 

\begin{tikzpicture}[x=0.55pt,y=0.55pt,yscale=-1,xscale=1]

\draw  [draw opacity=0] (382.33,208.44) .. controls (383.36,195.59) and (385.46,182.42) .. (388.73,169.14) .. controls (397.91,131.79) and (414.8,99.45) .. (435.74,75.48) -- (497.72,195.93) -- cycle ; \draw   (382.33,208.44) .. controls (383.36,195.59) and (385.46,182.42) .. (388.73,169.14) .. controls (397.91,131.79) and (414.8,99.45) .. (435.74,75.48) ;  
\draw  [draw opacity=0] (435.74,75.48) .. controls (447.28,81.23) and (458.71,88.1) .. (469.81,96.09) .. controls (501.03,118.56) and (524.72,146.31) .. (539.14,174.68) -- (404.25,187.19) -- cycle ; \draw   (435.74,75.48) .. controls (447.28,81.23) and (458.71,88.1) .. (469.81,96.09) .. controls (501.03,118.56) and (524.72,146.31) .. (539.14,174.68) ;  
\draw   (426,93.5) .. controls (426,86.04) and (432.04,80) .. (439.5,80) .. controls (446.96,80) and (453,86.04) .. (453,93.5) .. controls (453,100.96) and (446.96,107) .. (439.5,107) .. controls (432.04,107) and (426,100.96) .. (426,93.5) -- cycle ;
\draw  [draw opacity=0] (91.97,208.32) .. controls (93,195.47) and (95.1,182.29) .. (98.37,169.01) .. controls (107.55,131.66) and (124.44,99.33) .. (145.38,75.36) -- (207.36,195.8) -- cycle ; \draw   (91.97,208.32) .. controls (93,195.47) and (95.1,182.29) .. (98.37,169.01) .. controls (107.55,131.66) and (124.44,99.33) .. (145.38,75.36) ;  
\draw  [draw opacity=0] (145.38,75.36) .. controls (156.92,81.11) and (168.35,87.98) .. (179.45,95.97) .. controls (210.67,118.44) and (234.36,146.18) .. (248.78,174.56) -- (113.89,187.07) -- cycle ; \draw   (145.38,75.36) .. controls (156.92,81.11) and (168.35,87.98) .. (179.45,95.97) .. controls (210.67,118.44) and (234.36,146.18) .. (248.78,174.56) ;  
\draw    (253.33,141) -- (378.33,141) ;
\draw [shift={(380.33,141)}, rotate = 180] [color={rgb, 255:red, 0; green, 0; blue, 0 }  ][line width=0.75]    (10.93,-3.29) .. controls (6.95,-1.4) and (3.31,-0.3) .. (0,0) .. controls (3.31,0.3) and (6.95,1.4) .. (10.93,3.29)   ;
\draw  [draw opacity=0][line width=1.5]  (382.33,208.44) .. controls (383.36,195.59) and (385.46,182.42) .. (388.73,169.14) .. controls (396.69,136.75) and (410.45,108.13) .. (427.62,85.47) -- (497.72,195.93) -- cycle ; \draw  [color={rgb, 255:red, 208; green, 2; blue, 27 }  ,draw opacity=1 ][line width=1.5]  (382.33,208.44) .. controls (383.36,195.59) and (385.46,182.42) .. (388.73,169.14) .. controls (396.69,136.75) and (410.45,108.13) .. (427.62,85.47) ;  
\draw  [draw opacity=0][line width=1.5]  (448.59,82.38) .. controls (455.76,86.51) and (462.85,91.08) .. (469.81,96.09) .. controls (501.03,118.56) and (524.72,146.31) .. (539.14,174.68) -- (404.25,187.19) -- cycle ; \draw  [color={rgb, 255:red, 208; green, 2; blue, 27 }  ,draw opacity=1 ][line width=1.5]  (448.59,82.38) .. controls (455.76,86.51) and (462.85,91.08) .. (469.81,96.09) .. controls (501.03,118.56) and (524.72,146.31) .. (539.14,174.68) ;  
\draw  [draw opacity=0][line width=1.5]  (427.77,86.82) .. controls (430.26,82.45) and (435.09,79.66) .. (440.44,80.03) .. controls (443.44,80.24) and (446.15,81.42) .. (448.27,83.24) -- (439.5,93.5) -- cycle ; \draw  [color={rgb, 255:red, 208; green, 2; blue, 27 }  ,draw opacity=1 ][line width=1.5]  (427.77,86.82) .. controls (430.26,82.45) and (435.09,79.66) .. (440.44,80.03) .. controls (443.44,80.24) and (446.15,81.42) .. (448.27,83.24) ;  
\end{tikzpicture}
\caption{Representation of the insertion of a small enough disk to fix the non-differentiability of a point of the curved polygon approximation.}
\end{figure}
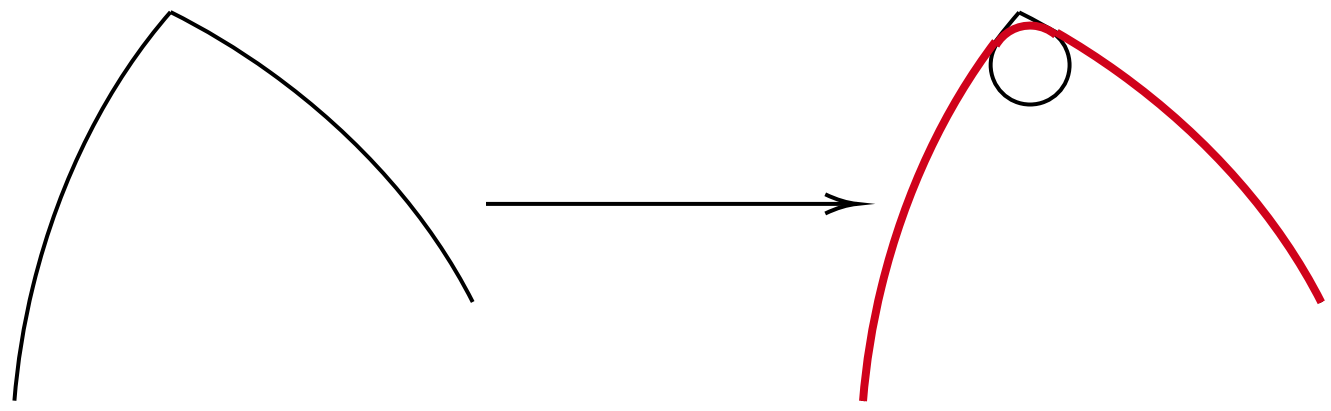

The only problem that could arise is if our original curve already had a corner, which could lead the approximating curved polygon to share that same corner, making the insertion of a small disk satisfying the conditions described above impossible. To overcome this particular case, we can instead approximate the original curve by slightly more inflated curved polygons, namely, ($\{\gamma_n+\delta B\}_n$). As $\delta\to 0$ and $n\to\infty$, these curves converge to the original curve ($\gamma$). However, the inflation ensures that the insertion of a disk with sufficiently small radius becomes possible while still containing ($\gamma$). Furthermore, this modification preserve perimeter relations between the curved polygons, i.e., if $L(\gamma_n)>L(\gamma_m)\implies L(\gamma_n+\delta B)>L(\gamma_m+\delta B)$.

The resulting arc remains arbitrarily close to the curved polygon. Furthermore, since the radius of the inserted circles satisfy $\varepsilon\ll1/\lambda$, the new arcs possess a curvature of at least $\lambda$, which means they preserve the $\lambda$-convexity of the overall curve while continuing to enclose $\gamma$. Hence, without loss of generality, we may assume that our approximated curves $\gamma_n$ are of class $\mathcal{C}^1$, implying that their tangents are defined at every point.

\section{Proofs of the main results}

The purpose of this section is to prove  Theorems \ref{theorem A} and \ref{theorem B}. To that end, we establish two complementary results that will allow us to assume additional regularity properties on arbitrary $\lambda$-convex standing arcs. 

\subsection{Convergence of the sequence of 1-convex lenses $\gamma^{*}_n$}

We begin by introducing a result that will allow us to simplify the proofs of Theorems \ref{theorem A} and \ref{theorem B} using the approximations described in Sections \ref{subsubsection: curved polygons} and \ref{subsubsection: small circles approxiamtion}.

\begin{proposition}
\label{proposition about convergence of lenses}
Let $\gamma$ be a $\lambda$-convex standing arc of finite length $L(\gamma)$, and let $\{\gamma_n\}$ be the sequence of curved polygons obtained by the approximation process described in section \ref{subsubsection: curved polygons}, such that $\gamma_n\xrightarrow{H}\gamma$ and $L(\gamma_n)\rightarrow L(\gamma)$.

For each $n$, let $\gamma_n^*$ be the standing $\lambda$-convex lens (respectively, the standing circular arc of radius $1/\lambda$) with $L(\gamma_n^*)=L(\gamma_n)$. Then $\gamma_n^*\xrightarrow{H}\gamma^*$, where $\gamma^*$ is the standing $\lambda$-convex lens (standing circular arc of radius 1/$\lambda$) with $L(\gamma^*)=L(\gamma)$.
\end{proposition}
This will allow us to solve the problem for the approximated $C^1$-smooth curves and extend that result to the limit of its sequence; which is of course our initial curve $\gamma$.

We prove Proposition \ref{proposition about convergence of lenses} for the case of the convergence of standing 1-convex lenses, since the proof for the convergence of standing circular arcs of radius 1 is similar and, in fact, simpler.

\begin{proof}[Proof of Proposition \ref{proposition about convergence of lenses}]
As will be recurrent, we prove the statement for $\lambda=1$, to simplify notation, and the general case $\lambda>0$ is obtained by scaling. 

Up to a horizontal translation we may assume that all lenses have the same contact point with the $x$-axis. Since all these lenses have the same curvature and are placed with the same contact point and the same axis of symmetry, the Hausdorff distance between them is the length of one of the segments shown in Figure \ref{Hausdorff distance between lenses}; depending on the lengths of the two lenses.
\\
\tikzset{every picture/.style={line width=0.75pt}} 

\begin{figure}[h]
\centering
\tikzset{every picture/.style={line width=0.75pt}} 

\begin{tikzpicture}[x=0.75pt,y=0.75pt,yscale=-1.5,xscale=1.5]

\draw  [draw opacity=0] (277.91,220.82) .. controls (260.45,239.82) and (223.17,252.89) .. (179.89,252.82) .. controls (136.39,252.74) and (98.72,239.41) .. (81.02,220.19) -- (179.37,195.95) -- cycle ; \draw  [color={rgb, 255:red, 208; green, 2; blue, 27 }  ,draw opacity=1 ] (277.91,220.82) .. controls (260.45,239.82) and (223.17,252.89) .. (179.89,252.82) .. controls (136.39,252.74) and (98.72,239.41) .. (81.02,220.19) ;  
\draw  [draw opacity=0] (81.02,220.19) .. controls (98.55,201.24) and (135.87,188.3) .. (179.15,188.53) .. controls (222,188.76) and (259.15,201.82) .. (277.11,220.64) -- (179.47,245.4) -- cycle ; \draw  [color={rgb, 255:red, 208; green, 2; blue, 27 }  ,draw opacity=1 ] (81.02,220.19) .. controls (98.55,201.24) and (135.87,188.3) .. (179.15,188.53) .. controls (222,188.76) and (259.15,201.82) .. (277.11,220.64) ;  
\draw  [draw opacity=0] (92.31,230.72) .. controls (112.16,217.14) and (143.65,208.47) .. (179.15,208.66) .. controls (214.15,208.84) and (245.35,217.59) .. (265.39,231.05) -- (179.47,265.52) -- cycle ; \draw  [color={rgb, 255:red, 74; green, 144; blue, 226 }  ,draw opacity=1 ] (92.31,230.72) .. controls (112.16,217.14) and (143.65,208.47) .. (179.15,208.66) .. controls (214.15,208.84) and (245.35,217.59) .. (265.39,231.05) ;  
\draw    (180,188.6) -- (180,208.6) ;
\draw    (27,252.6) -- (329,253.6) ;

\draw (181.99,191.6) node [anchor=north west][inner sep=0.75pt]  [font=\small,rotate=-0.11] [align=left] {$\displaystyle d_{H,1}$};

\draw (265.39,231.05) -- (277.91,220.82);

\draw (271.65,225.94) node [anchor=north west] [inner sep=2pt]
[font=\small] {$\displaystyle d_{H,2}$};

\end{tikzpicture}

\caption{Illustration of the possible Hausdorff distances between two $\lambda$-convex lens.}
\label{Hausdorff distance between lenses}
\end{figure}
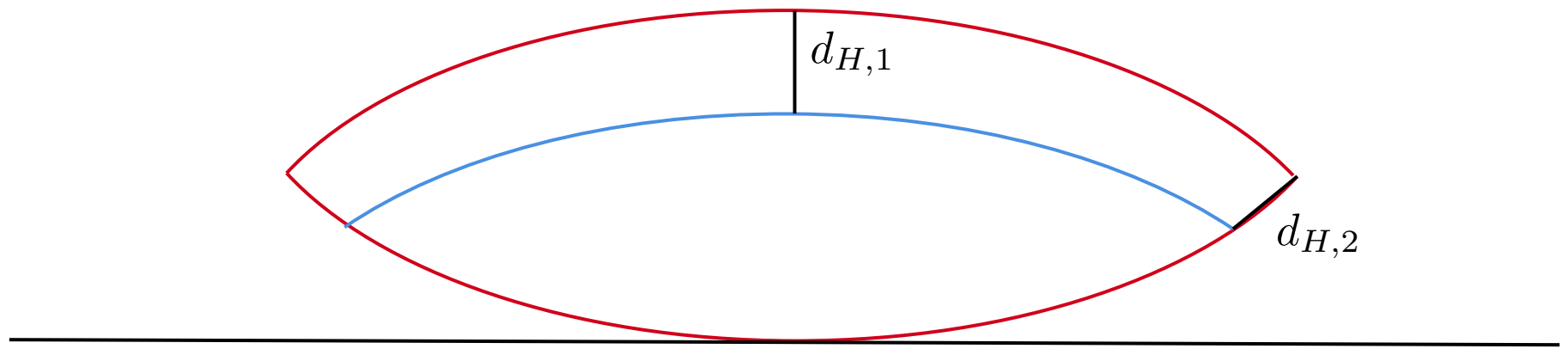
It is clear that the sequence $\{L(\gamma_n^*)\}$ converges to the length of our original curve $\gamma$. By hypothesis we have that $\gamma_n\xrightarrow{H}\gamma$ and $L(\gamma_n)\rightarrow L(\gamma)$. By construction we have $L(\gamma_n)=L(\gamma_n^*)$ $\forall n$, so we can trivially rewrite the last convergence as $L(\gamma_n^*)\rightarrow L(\gamma^*)=L(\gamma)$. We can now ensure that $\{L(\gamma_n^*)\}$ is a convergent sequence and, as such, $|L(\gamma_n^*)-L(\gamma)|\xrightarrow{n\rightarrow\infty}$ 0. From here it is simple to see what we want. We will compute the two possible Hausdorff distances shown in Figure \ref{Hausdorff distance between lenses} and show they tend to 0 as $n\rightarrow\infty$.

Let $L_n=L(\gamma_n)=L(\gamma_n^*)$. We compute the highest point of the 1-convex lens of length $L_n$. Such a lens consists of the symmetric intersection of two circular arcs of curvature 1, each of length $L_n/2$. Let $\alpha$ be the angle corresponding to half of one of these arcs (see Figure \ref{drawin for calculation} (a)).
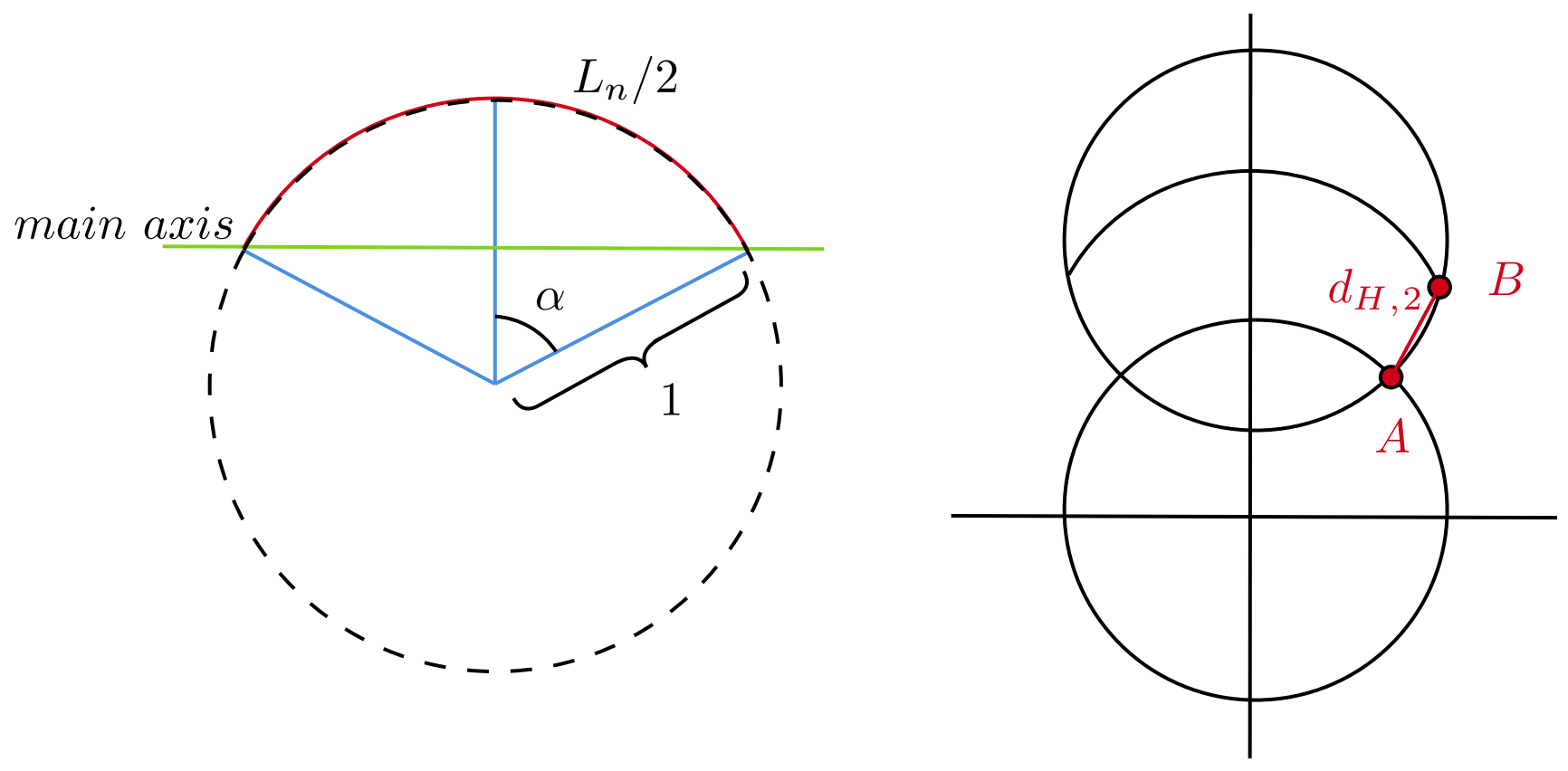
\begin{figure}[h!]
\centering

\tikzset{every picture/.style={line width=0.75pt}} 

\begin{tikzpicture}[x=0.55pt,y=0.55pt,yscale=-1,xscale=1]

\draw [color={rgb, 255:red, 74; green, 144; blue, 226 }  ,draw opacity=1 ]   (204,60) -- (204,168) ;
\draw [color={rgb, 255:red, 74; green, 144; blue, 226 }  ,draw opacity=1 ]   (299.46,118.35) -- (204,168) ;
\draw [color={rgb, 255:red, 74; green, 144; blue, 226 }  ,draw opacity=1 ]   (109.04,117.52) -- (204,168) ;
\draw [color={rgb, 255:red, 126; green, 211; blue, 33 }  ,draw opacity=1 ]   (78.4,116) -- (328.4,117) ;
\draw   (211,173.4) .. controls (213.25,177.49) and (216.42,178.4) .. (220.51,176.15) -- (249.12,160.36) .. controls (254.96,157.14) and (259.01,157.57) .. (261.26,161.66) .. controls (259.01,157.57) and (260.8,153.92) .. (266.64,150.69)(264.01,152.14) -- (295.25,134.91) .. controls (299.34,132.65) and (300.25,129.48) .. (298,125.39) ;
\draw  [draw opacity=0] (108.88,116.81) .. controls (127.12,82.98) and (162.88,60) .. (204,60) .. controls (245.44,60) and (281.43,83.34) .. (299.54,117.59) -- (204,168) -- cycle ; \draw  [color={rgb, 255:red, 208; green, 2; blue, 27 }  ,draw opacity=1 ] (108.88,116.81) .. controls (127.12,82.98) and (162.88,60) .. (204,60) .. controls (245.44,60) and (281.43,83.34) .. (299.54,117.59) ;  
\draw  [draw opacity=0] (204.02,142.46) .. controls (208.34,142.72) and (212.74,143.97) .. (216.86,146.3) .. controls (221.14,148.72) and (224.59,152.02) .. (227.05,155.77) -- (204.67,167.84) -- cycle ; \draw   (204.02,142.46) .. controls (208.34,142.72) and (212.74,143.97) .. (216.86,146.3) .. controls (221.14,148.72) and (224.59,152.02) .. (227.05,155.77) ;  
\draw  [draw opacity=0][dash pattern={on 4.5pt off 4.5pt}] (109.04,117.52) .. controls (127.28,83.7) and (163.03,60.72) .. (204.16,60.72) .. controls (263.8,60.72) and (312.16,109.07) .. (312.16,168.72) .. controls (312.16,228.36) and (263.8,276.72) .. (204.16,276.72) .. controls (144.51,276.72) and (96.16,228.36) .. (96.16,168.72) .. controls (96.16,150.26) and (100.79,132.88) .. (108.95,117.69) -- (204.16,168.72) -- cycle ; \draw  [color={rgb, 255:red, 0; green, 0; blue, 0 }  ,draw opacity=1 ][dash pattern={on 4.5pt off 4.5pt}] (109.04,117.52) .. controls (127.28,83.7) and (163.03,60.72) .. (204.16,60.72) .. controls (263.8,60.72) and (312.16,109.07) .. (312.16,168.72) .. controls (312.16,228.36) and (263.8,276.72) .. (204.16,276.72) .. controls (144.51,276.72) and (96.16,228.36) .. (96.16,168.72) .. controls (96.16,150.26) and (100.79,132.88) .. (108.95,117.69) ;  
\draw   (419.17,113.69) .. controls (419.17,74.01) and (451.56,41.84) .. (491.51,41.84) .. controls (531.46,41.84) and (563.85,74.01) .. (563.85,113.69) .. controls (563.85,153.38) and (531.46,185.55) .. (491.51,185.55) .. controls (451.56,185.55) and (419.17,153.38) .. (419.17,113.69) -- cycle ;
\draw   (419.17,215.65) .. controls (419.17,175.97) and (451.56,143.8) .. (491.51,143.8) .. controls (531.46,143.8) and (563.85,175.97) .. (563.85,215.65) .. controls (563.85,255.34) and (531.46,287.51) .. (491.51,287.51) .. controls (451.56,287.51) and (419.17,255.34) .. (419.17,215.65) -- cycle ;
\draw    (489.56,28) -- (489.31,309.6) ;
\draw    (376.4,217.84) -- (605.4,218.57) ;
\draw  [draw opacity=0] (420.9,126.72) .. controls (433.77,103.65) and (459.32,87.8) .. (488.9,87.47) .. controls (521.02,87.12) and (548.8,105.18) .. (560.92,131.25) -- (489.7,160.69) -- cycle ; \draw   (420.9,126.72) .. controls (433.77,103.65) and (459.32,87.8) .. (488.9,87.47) .. controls (521.02,87.12) and (548.8,105.18) .. (560.92,131.25) ;  
\draw  [fill={rgb, 255:red, 208; green, 2; blue, 27 }  ,fill opacity=1 ] (556.96,131.42) .. controls (556.96,129.21) and (558.77,127.41) .. (561,127.41) .. controls (563.22,127.41) and (565.03,129.21) .. (565.03,131.42) .. controls (565.03,133.63) and (563.22,135.42) .. (561,135.42) .. controls (558.77,135.42) and (556.96,133.63) .. (556.96,131.42) -- cycle ;
\draw  [fill={rgb, 255:red, 208; green, 2; blue, 27 }  ,fill opacity=1 ] (538.64,165.4) .. controls (538.64,163.19) and (540.44,161.4) .. (542.67,161.4) .. controls (544.89,161.4) and (546.7,163.19) .. (546.7,165.4) .. controls (546.7,167.62) and (544.89,169.41) .. (542.67,169.41) .. controls (540.44,169.41) and (538.64,167.62) .. (538.64,165.4) -- cycle ;
\draw [color={rgb, 255:red, 208; green, 2; blue, 27 }  ,draw opacity=1 ]   (561,131.42) -- (542.67,165.4) ;

\draw (264,166) node [anchor=north west][inner sep=0.75pt]  [font=\small] [align=left] {$\displaystyle 1$};
\draw (217,130) node [anchor=north west][inner sep=0.75pt]  [font=\small] [align=left] {$\displaystyle \alpha $};
\draw (231,42) node [anchor=north west][inner sep=0.75pt]  [font=\small] [align=left] {$\displaystyle L_{n} /2$};
\draw (29,86.4) node [anchor=north west][inner sep=0.75pt]    {$$};
\draw (20,99) node [anchor=north west][inner sep=0.75pt]  [font=\small] [align=left] {$\displaystyle main\ axis$};
\draw (576.8,119.91) node [anchor=north west][inner sep=0.75pt]  [font=\small,color={rgb, 255:red, 208; green, 2; blue, 27 }  ,opacity=1 ] [align=left] {$\displaystyle B$};
\draw (534.26,179.39) node [anchor=north west][inner sep=0.75pt]  [font=\small,color={rgb, 255:red, 208; green, 2; blue, 27 }  ,opacity=1 ] [align=left] {$\displaystyle A$};
\draw (516.59,122.44) node [anchor=north west][inner sep=0.75pt]  [font=\small,color={rgb, 255:red, 208; green, 2; blue, 27 }  ,opacity=1 ] [align=left] {$\displaystyle d_{H}{}_{,}{}_{2}$};
\draw (192,325) node [anchor=north west][inner sep=0.75pt]   [align=left] {$\displaystyle ( a)$};
\draw (481,325) node [anchor=north west][inner sep=0.75pt]   [align=left] {$\displaystyle ( b)$};

\end{tikzpicture}

\caption{(a) Illustration used in the calculation of the height of a 1-convex lens. We have drawn the full circle (dotted line) of radius 1 that corresponds to the upper-half of the 1-convex lens (red line). The green line corresponds to the main axis of the 1-convex lens. (b) Illustration for the computation of the distance $d_{H,2}$ shown in Figure \ref{Hausdorff distance between lenses}. The intersection of these two circumferences of radius 1 represent the smaller of the 1-convex lenses, so $A$ is its right most corner, while $B$ is the right most corner of the larger lens.}
\label{drawin for calculation}
\end{figure}

Since the radius $r$ is 1 we have:
\begin{equation*}
    \frac{L_n}{4}=r\cdot\alpha=1\cdot\alpha\implies\alpha=\frac{L_n}{4},
\end{equation*}
and by some basic trigonometry we conclude that the distance between the main axis and the highest point of the 1-convex lens is $1-\text{cos}(L_n/4)$. Due to vertical symmetry, the 1-convex lens has with respect to its main axis the height of the lens will be twice the distance between its highest point and its main axis; meaning we obtain a total height of the lens of $2-2\text{cos}(L_n/4)$. 

Due to the assumed initial normalization, $d_{H,1}$ between two such lenses equals the difference of their heights:
\begin{equation*}
    d_{H,1}(\gamma_n^*,\gamma^*)=2-2+2\left(\text{cos}\left(\frac{L(\gamma)}{4}\right)-\text{cos}\left(\frac{L_n}{4}\right)\right)=2\left(\text{cos}\left(\frac{L(\gamma)}{4}\right)-\text{cos}\left(\frac{L_n}{4}\right)\right).
\end{equation*}

Therefore, since $L_n\rightarrow L(\gamma)$ and cosine is a continuous function, it is clear that the value of the last expression approaches 0 when $n$ grows.

Next, we compute $d_{H,2}$. We will use Figure \ref{drawin for calculation} (b). Observe that the point $A$ in the said Figure can be seen as the intersection of two circumferences of radius 1, one centered at $(0,0)$ and the other at $(0.r)$, where $r$ is 2 minus the width of the smaller lens, which we already know from the previous computation. Therefore, we know that $r=2\cos\left(\frac{L(\gamma)}{4}\right)$. From here it is quite simple to find the exact coordinates of $A$ and $B$ by using the equations of the circles and the width of the lenses. We obtain:
\begin{gather*}
    A=\left(\sin\left(\frac{L(\gamma)}{4}\right), \cos\left(\frac{L(\gamma)}{4}\right)\right), \\
    B=\left(\sin\left(\frac{L_n}{4}\right), 2\cos\left(\frac{L(\gamma)}{4}\right)-\cos\left(\frac{L_n}{4}\right)\right).
\end{gather*}
Therefore, the expression of the desired distance is:
\begin{equation*}
    d_{H,2}(\gamma_n^*,\gamma^*)=\sqrt{\left(\sin\frac{L(\gamma)}{4}-\sin\frac{L_n}{4}\right)^2+\left(\cos\frac{L_n}{4}-\cos\frac{L(\gamma)}{4}\right)^2}.
\end{equation*}
Therefore, since $L_n\rightarrow L(\gamma)$ and cosine and sine are a continuous functions, the last expression approaches 0 when $n$ grows.

As a result, we can conclude that $d_{H,1}(\gamma_n^*,\gamma^*), d_{H,2}(\gamma_n^*,\gamma^*)\xrightarrow{n\rightarrow\infty}$ 0; which means, by definition of the Hausdorff metric, that the sequence of 1-convex lenses converges to $\gamma^*$.
\end{proof}
The importance of this step of the proof revolves around the fact that the center of mass of these curves is continuous with respect to the weak convergence of their corresponding arc-length measures (see \cite{AmbrosioFuscoPallara2000}). Since $\gamma_n \xrightarrow{H} \gamma$ and $\gamma_n^* \xrightarrow{H} \gamma^*$, and the lengths of the corresponding curves converge, we obtain:
\begin{gather*}
CM(\gamma_n)\rightarrow CM(\gamma), \\
CM(\gamma_n^*)\rightarrow CM(\gamma^*).
\end{gather*}
If for each $n$ we have $CM(\gamma_n)\ge CM(\gamma_n^*)$ passing to the limit we get:
\begin{equation*}
    CM(\gamma)\ge CM(\gamma^*).
\end{equation*}

So, once Theorem \ref{theorem A} or Theorem \ref{theorem B} are proven for $\gamma_n$, the  results follow for our original curve $\gamma$. 

\subsection{Proof of Theorem \ref{theorem B}}

\begin{proof}[Proof of Theorem \ref{theorem B}]

By scaling, we may assume $\lambda=1$. Using the curved polygon approximation and the small circle approximation from Sections \ref{subsubsection: curved polygons} and \ref{subsubsection: small circles approxiamtion} and Proposition \ref{proposition about convergence of lenses}, it is enough to prove the result for $\mathcal{C}^1$ standing arcs whose tangent vector is defined everywhere.

Let $L=L(\gamma)$, and let $\tilde{\gamma}$ be the standing circular arc of radius $1$ and length $L$.

Parametrize $\gamma$ by arc length from its left endpoint to its right endpoint so 
\[
\gamma(s)=(x(s),y(s)) \quad \text{for} \quad s\in[0,L].
\]
Let \(\alpha(s)\) be the oriented angle of the tangent vector with the positive \(x\)-axis, so that
\[
\gamma'(s)=(\cos\alpha(s),\sin\alpha(s)).
\]
Since \(\gamma\) is convex and is traversed from left to right, \(\alpha(s)\) is non-increasing, and with
this orientation we have that its curvature is
\[
k(s)=-\alpha'(s).
\]
In particular, since \(\gamma\) is \(1\)-convex, \(k(s)\ge1\), and therefore $\alpha'(s)$ is upper bounded by -1.

Observe that this tangent angle $\alpha(s)$ is related to the contact angles (Definition \ref{contact angles}) as $\alpha(0)=\theta_1(\gamma)$ and $\alpha(L)=-\theta_2(\gamma)$. Since \(\theta_1(\gamma),\theta_2(\gamma)\le \pi/2\) under the conditions of Theorem \ref{theorem B}, we have $
\alpha(s)\in[-\pi/2,\pi/2]
$
for every \(s\in[0,L]\). Moreover, using $k(s)\ge1$ we can upper bound the length of $\gamma$ as
\[
L\le \int_0^L k(s)\,ds
=
\alpha(0)-\alpha(L)
=
\theta_1(\gamma)+\theta_2(\gamma)
\le \pi.
\]

We can do the same for the standing circular arc of radius 1 and length $L$, $\tilde{\gamma}$. In this case the expression of the tangent angle is known explicitly:
\[
\widetilde{\alpha}(s)=\frac{L}{2}-s,
\]
which gives $\tilde{\alpha}'(s)=-1$.

Since \(L\le\pi\) and $s\in[0,L]$, we also have
\(
\widetilde{\alpha}(s)\in[-\pi/2,\pi/2].
\)

From here, we can define
\[
g(s):=\alpha(s)-\widetilde{\alpha}(s).
\]
Differentiating this last expression we get
\[
g'(s)=\alpha'(s)-\widetilde{\alpha}'(s)
=
\alpha'(s)+1
\le0,
\]
where for the last inequality we have used $\alpha'(s)\le-1$. Thus $g$ is non-increasing.

Now set
\[
f(s):=\sin\alpha(s)-\sin\widetilde{\alpha}(s).
\]
Since both \(\alpha(s)\) and \(\widetilde{\alpha}(s)\) lie in \([-\pi/2,\pi/2]\), and since the sine
function is increasing on this interval, the sign of \(f(s)\) is the same as the sign of
\(g(s)\) ($\alpha(s)>\tilde{\alpha}(s)\iff\sin\alpha(s)>\sin\tilde{\alpha}(s))$. Therefore \(f\) changes sign at most once, and if it does, it changes from positive
to negative.

Using $\gamma'(s)=(\cos\alpha(s),\sin\alpha(s))$ and $\tilde{\gamma}'(s)=(\cos\tilde{\alpha}(s),\sin\tilde{\alpha}(s))$, we can write
\begin{equation}
\label{thm b: expressions y}
y(s)=\int_0^s \sin\alpha(t)\,dt, \qquad
\tilde{y}(s)=\int_0^s \sin\tilde{\alpha}(t)\,dt.
\end{equation}
Now, since both $\gamma$ and $\tilde{\gamma}$ are standing arcs, they start and end on the $x$-axis. In particular, $y(L)=\tilde{y}(L)=0$. Evaluating the expressions above at $s=L$, we obtain
\[
\int_0^L \sin\alpha(s)\,ds=0,
\qquad
\int_0^L \sin\tilde{\alpha}(s)\,ds=0.
\]
Therefore, recalling the definition of $f(s)$, we get
\begin{equation}
\label{thm b: int f=0}
\int_0^L f(s)\,ds=0.
\end{equation}
Applying Fubini's theorem gives
\[
\int_0^L y(s)\,ds
=
\int_0^L (L-t)\sin\alpha(t)\,dt,
\qquad
\int_0^L \widetilde y(s)\,ds
=
\int_0^L (L-t)\sin\widetilde{\alpha}(t)\,dt,
\]
and subtracting, we obtain
\begin{equation}
\label{thm b ÚLTIMA}
\int_0^L \bigl(y(s)-\widetilde y(s)\bigr)\,ds
=
\int_0^L (L-t)f(t)\,dt.
\end{equation}
If \(f\) does not change sign from \eqref{thm b: int f=0} we get that $f\equiv0$, which implies that $\sin\alpha(s)=\sin\tilde\alpha(s)$ $\forall s\in[0,L]$. Therefore, from \eqref{thm b: expressions y} we have $y(s)=\tilde{y}(s)$ $\forall s\in[0,L]$ and the centers of mass of $\gamma$ and $\tilde\gamma$ are at the same height. As a matter of fact, in this case we would have $\gamma=\tilde\gamma$.

Now suppose $f$ does change sign. Let \(a\in[0,L]\) be such that \(f\ge0\) on \([0,a]\) and \(f\le0\) on \([a,L]\). Since \(L-t\) is decreasing, we have
\[
L-t\ge L-a \quad \text{for } \ t\in[0,a],
\]
and
\[
L-t\le L-a \quad \text{for } \ t\in[a,L].
\]
Using also \(f\ge0\) on \([0,a]\), \(f\le0\) on \([a,L]\), and \eqref{thm b: int f=0}, we get
\begin{align*}
\int_0^L (L-t)f(t)\,dt
&=
\int_0^a (L-t)f(t)\,dt
+
\int_a^L (L-t)f(t)\,dt\\
&\ge
(L-a)\int_0^a f(t)\,dt
+
(L-a)\int_a^L f(t)\,dt\\
&=
(L-a)\int_0^L f(t)\,dt=0.
\end{align*}
Going back to \eqref{thm b ÚLTIMA} we have that
\[
\int^L_0(y(s)-\tilde{y}(s))ds=\int^L_0(L-t)f(t)dt\ge0\implies\int^L_0y(s)ds\ge\int^L_0\tilde{y}(s)ds.
\]
Dividing by \(L\) and recalling the equation for the center of mass \eqref{center of mass eqn}, we obtain
\[
CM(\gamma)\ge CM(\widetilde{\gamma}).
\]
This proves Theorem \ref{theorem B} for \(\lambda=1\). The general case, as mentioned, follows by scaling.

\end{proof}

\subsection{Reduction to an equivalent case with vertical symmetry}

From our approximated arc, which we will call $\gamma$ for the rest of the paper to simplify notation, we can describe the half perimeter vertical line $HVP(\gamma)$ as defined in \ref{def half vertical}. Our objective now is to reduce our problem to one with some sort of vertical symmetry; which will simplify the proof of Theorem \ref{theorem A}. For this, we will prove the following statement:

\begin{proposition}\label{symmetry proposition}
Let $\gamma\subset(\mathbb R^2)^+$ be a $\lambda$--convex standing arc of finite length $L(\gamma)$
with constant linear density. 
Then there exists a $\lambda$--convex standing arc $\gamma^r$ such that
\begin{enumerate}
  \item $\gamma^r$ is symmetric with respect to $HVP(\gamma)$,
  \item $L(\gamma^r)=L(\gamma)$,
  \item $CM(\gamma^r)\le CM(\gamma)$.
\end{enumerate}
In particular, $HVP(\gamma^r)=HVP(\gamma)$.
\end{proposition}

\begin{proof}

By the curved--polygon approximation and the small-circle regularization used throughout the article,
we may assume, without loss of generality, that $\gamma$ is $C^1$.

Let $HVP(\gamma)$ be its half--perimeter vertical line and write
$\gamma=\gamma_-\cup\gamma_+$ for the decomposition into the left and right halves determined by $HVP(\gamma)$, respectively,
so that $L(\gamma_-)=L(\gamma_+)=\tfrac12L(\gamma)$.

 Let $p$ be the highest point of $\gamma\cap HVP(\gamma)$.
Let $\alpha_-$ and $\alpha_+$ be the left and right \emph{section angles at $p$},
measured with the convention of Definition \ref{def section angles} (see Figure \ref{fig:section angles}): they are the angles between the positive $y$-axis
and the tangent line at $p$ when approaching from $\gamma_-$ and $\gamma_+$, respectively.

Since $HVP(\gamma)$ splits $\gamma$ into two arcs of equal length,
\begin{equation}\label{eq:CM-half-average}
CM(\gamma)=\frac{CM(\gamma_-)+CM(\gamma_+)}{2}.
\end{equation}

We distinguish cases according to whether the section angles at $p$ agree.

\medskip
\begin{enumerate}
\item $\alpha_-=\alpha_+$.
Assume, without loss of generality, $CM(\gamma_-)\le CM(\gamma_+)$.
Define $\gamma^r$ by keeping $\gamma_-$ and replacing $\gamma_+$ by the reflection of $\gamma_-$
with respect to $HVP(\gamma)$.
Reflection is an isometry, hence preserves length and $\lambda$--convexity on each piece.
Because $\alpha_-=\alpha_+$, the tangents match at $p$, so no corner is created there.
Therefore $\gamma^r$ is $\lambda$--convex and standing, $L(\gamma^r)=L(\gamma)$, and
$$CM(\gamma^r)=\frac{CM(\gamma_-)+CM(\gamma_-)}{2}=CM(\gamma_-)\le CM(\gamma),$$
by \eqref{eq:CM-half-average}. 
\medskip
\item $\alpha_-\neq\alpha_+$.
Without loss of generality assume $\alpha_-<\alpha_+$.
We further split into two subcases according to which half has smaller center of mass.

\medskip
\begin{enumerate}
\item $CM(\gamma_-)\le CM(\gamma_+)$.
As in case 1, we replace $\gamma_+$ by the reflection of $\gamma_-$
with respect to $HVP(\gamma)$.
Reflection preserves arc length and $\lambda$--convexity away from the gluing point.
Since $\alpha_-<\alpha_+$, after reflection both halves meet at $p$
with the same section angle $\alpha_-$ with respect to the $y$--axis.
Hence no inward corner is created (at worst an outward corner may appear,
which is compatible with $\lambda$--convexity).

We thus obtain a $\lambda$--convex standing arc $\gamma^r$ symmetric with respect to $HVP(\gamma)$ with
$L(\gamma^r)=L(\gamma)$ and, as above,
$$CM(\gamma^r)=CM(\gamma_-)\le CM(\gamma).$$

\medskip
\item $CM(\gamma_-)> CM(\gamma_+)$.

In this case we cannot reflect the entire right half $\gamma_+$
across $HVP(\gamma)$, since doing so would typically create
an inward corner at $p$ and thus violate $\lambda$--convexity.
Instead, we reflect only a suitable lower portion of the curve,
leaving the upper part unchanged. The goal is to determine a height
at which this partial symmetrization preserves $\lambda$--convexity.
We determine the appropriate gluing height as follows.

\medskip

Let $\ell:=L(\gamma)/2$ and parametrize each half $\gamma_\pm$
by arc length starting at the highest point $p$
and moving \emph{downwards} toward the $x$--axis:
\[
\gamma_\pm(s)=(x_\pm(s),y_\pm(s)), \qquad s\in[0,\ell],
\]
with $\gamma_\pm(0)=p$ and $y_\pm(\ell)=0$.

For each $s$, let $\varphi_\pm(s)\in[0,\pi]$ denote the angle between
the positive $y$--axis and the tangent to $\gamma_\pm$ at $\gamma_\pm(s)$.
Since the curve is parametrized from top to bottom, its vertical component
is negative, so the unit tangent vector can be written as
$\displaystyle T_\pm(s)=\bigl(\sin\varphi_\pm(s),\, -\cos\varphi_\pm(s)\bigr).$
In particular,
\begin{equation}\label{eq:yprime-cos}
y_\pm'(s)=-\cos\bigl(\varphi_\pm(s)\bigr).
\end{equation}

Define the height difference $H(s):=y_+(s)-y_-(s).$ Then $H$ is continuous on $[0,\ell]$ and satisfies $H(0)=0=H(\ell)$.

Let
\[
U:=\{s\in(0,\ell): H(s)>0\}.
\]
Since $\alpha_-<\alpha_+$ and the cosine function is strictly decreasing on $[0,\pi]$,
\[
H'(0)=y_+'(0)-y_-'(0)=\cos\alpha_- - \cos\alpha_+>0,
\]
It follows from the definition of the derivative that $H(s)>0$ for all sufficiently small $s>0$. Thus $U\neq\emptyset$. Since $U\subset(0,\ell)$, it is bounded above,
and therefore it admits a supremum. Set
\begin{equation}\label{eq:def-s0}
s_0:=\sup U \in (0,\ell].
\end{equation}
Since no point $s>s_0$ can belong to $U$,  $H(s)>0$ for all $s\in(0,s_0)$ and $H(s)\le 0$ for all $s\in(s_0,\ell]$.
By continuity we then have
\begin{equation}\label{eq:H-s0-zero}
H(s_0)=y_+(s_0)-y_-(s_0)=0.
\end{equation}
We perform the reflection at level $s_0$:
we keep $\gamma_\pm|_{[0,s_0]}$ unchanged and replace the lower left part
$\gamma_-|_{[s_0,\ell]}$ by the reflection of
$\gamma_+|_{[s_0,\ell]}$ with respect to $HVP(\gamma)$.
The two pieces meet at $\gamma_\pm(s_0)$,
which lie at the same height by \eqref{eq:H-s0-zero}.

At this gluing point, the incoming tangent from the reflected lower piece
has angle $\varphi_+(s_0)$, while the outgoing tangent of the
unchanged upper left piece has angle $\varphi_-(s_0)$.
To avoid creating an inward corner, and hence to preserve
$\lambda$-convexity, it is necessary and sufficient that
\begin{equation}\label{eq:good-angle}
\varphi_-(s_0)\ge \varphi_+(s_0).
\end{equation}

We now verify that this inequality holds.
Since $H(s_0)=0$ and $H(s)>0$ for $s<s_0$ while $H(s)\le 0$ for $s>s_0$,
$H'(s_0)$ cannot be positive. Using $y_\pm'(s)=-\cos(\varphi_\pm(s))$,
$$
H'(s_0)=y_+'(s_0)-y_-'(s_0)
=
\cos\varphi_-(s_0)-\cos\varphi_+(s_0)\leq 0
$$
implies that
\[
\cos\varphi_-(s_0)\le \cos\varphi_+(s_0),
\]
and since the cosine function is strictly decreasing on $[0,\pi]$,
\[
\varphi_-(s_0)\ge \varphi_+(s_0),
\]
i.e condition \eqref{eq:good-angle} holds.
\medskip
Now we will perform a partial symmetrization below height $y_0:=y_-(s_0)=y_+(s_0)$ (equal height by \eqref{eq:H-s0-zero}).
Split each half into an upper and a lower part at $s_0$:
\[
\gamma_\pm^{\mathrm{up}}:=\gamma_\pm([0,s_0]),
\qquad
\gamma_\pm^{\mathrm{low}}:=\gamma_\pm([s_0,\ell]).
\]
Note that $L(\gamma_-^{\mathrm{low}})=L(\gamma_+^{\mathrm{low}})=\ell-s_0$.

Define a new arc $\widetilde\gamma$ as follows:
\begin{itemize}
\item keep the entire upper part $\gamma_-^{\mathrm{up}}\cup\gamma_+^{\mathrm{up}}$ unchanged;
\item replace the lower left part $\gamma_-^{\mathrm{low}}$ by the reflection of $\gamma_+^{\mathrm{low}}$
with respect to $HVP(\gamma)$ (and keep the lower right part as the reflection of that new left lower part).
\end{itemize}
By construction, $\widetilde\gamma$ is symmetric with respect to $HVP(\gamma)$ below height $y_0$,
and $L(\widetilde\gamma)=L(\gamma)$.

\smallskip
We check now that $\widetilde\gamma$ is $\lambda$-convex.
Reflection preserves $\lambda$-convexity away from the glue points.
At the glue points on the line $y=y_0$,
the incoming tangent from the reflected lower piece has the same angle with the $y$--axis
as the original right lower tangent, namely $\varphi_+(s_0)$, while the outgoing tangent of the unchanged
upper left piece has angle $\varphi_-(s_0)$.
The inequality $\varphi_-(s_0)\ge\varphi_+(s_0)$  prevents the creation of an inward corner, hence the supporting disk condition for $\lambda$-convexity
is preserved at the junction. Therefore $\widetilde\gamma$ is $\lambda$-convex and standing. 

\smallskip
We check now that the center of mass of $\widetilde\gamma$ does not increase.
By definition of $s_0$ we have $H(s)\le 0$ for all $s\in[s_0,\ell]$, i.e.
$\displaystyle
y_+(s)\le y_-(s)\;\forall s\in[s_0,\ell].
$
Integrating over $[s_0,\ell]$ yields
\[
\int_{s_0}^{\ell} y_+(s)\,ds \ \le\ \int_{s_0}^{\ell} y_-(s)\,ds,
\]
so the lower right piece has no larger average height than the lower left piece.
Since $\widetilde\gamma$ is obtained by replacing $\gamma_-^{\mathrm{low}}$ with a reflected copy of
$\gamma_+^{\mathrm{low}}$ of the same length, it follows that $CM(\widetilde\gamma)\le CM(\gamma)$.
\medskip

The curve $\widetilde\gamma$ is symmetric with respect to $HVP(\gamma)$
below height $y_0$. In particular, the lower left and lower right parts
are mirror images of each other and therefore contribute equally to the
center of mass.

Hence any difference between the centers of mass of the two halves of
$\widetilde\gamma$ is entirely determined by the upper parts.
Applying the argument of Case~2a to $\widetilde\gamma$,
we reflect the upper half with larger center of mass across $HVP(\gamma)$.
As in Case~2a, this reflection preserves $\lambda$-convexity.

We denote by $\gamma^r$ the resulting curve.
Then $\gamma^r$ is symmetric with respect to $HVP(\gamma)$,
$L(\gamma^r)=L(\widetilde\gamma)=L(\gamma)$,
and
\[
CM(\gamma^r)\le CM(\widetilde\gamma)\le CM(\gamma).
\]This concludes case 2b) and finishes the proof of the proposition.
\end{enumerate}
\end{enumerate}

\end{proof}
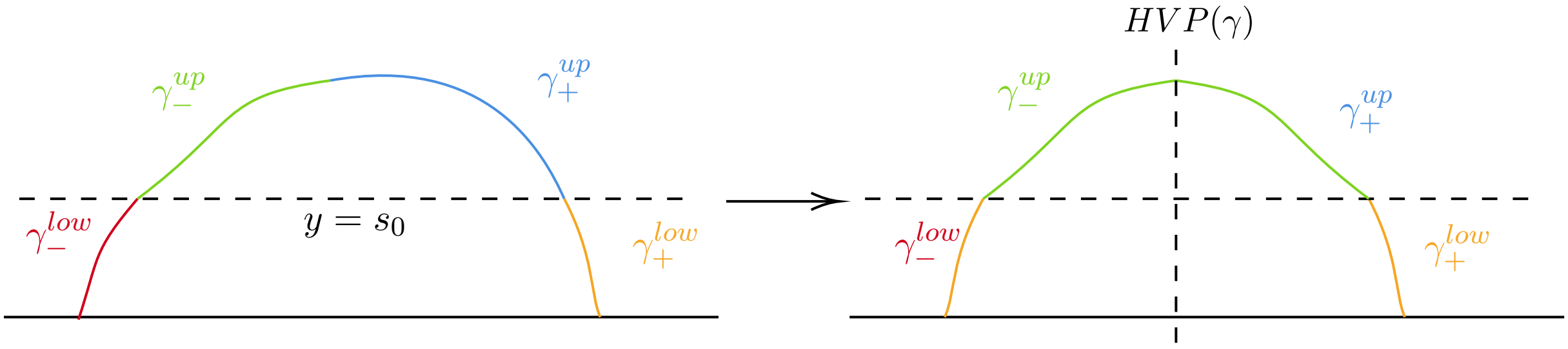
\begin{figure}[h]
\centering
\tikzset{every picture/.style={line width=0.75pt}} 

\begin{tikzpicture}[x=0.75pt,y=0.75pt,yscale=-1,xscale=1]

\draw    (42,161.6) -- (320,161.6) ;
\draw [color={rgb, 255:red, 208; green, 2; blue, 27 }  ,draw opacity=1 ]   (71,162.6) .. controls (79,138.6) and (76,136.6) .. (94,115.6) ;
\draw [color={rgb, 255:red, 126; green, 211; blue, 33 }  ,draw opacity=1 ]   (94,115.6) .. controls (134,85.6) and (125,75.6) .. (169,69.6) ;
\draw [color={rgb, 255:red, 74; green, 144; blue, 226 }  ,draw opacity=1 ]   (169,69.6) .. controls (230,58.6) and (252,97.6) .. (260,115.6) ;
\draw [color={rgb, 255:red, 245; green, 166; blue, 35 }  ,draw opacity=1 ]   (260,115.6) .. controls (272,138.6) and (270,152.6) .. (274,161.6) ;
\draw  [dash pattern={on 4.5pt off 4.5pt}]  (48,115.6) -- (307,115.6) ;
\draw    (323,116.6) -- (365,116.6) ;
\draw [shift={(367,116.6)}, rotate = 180] [color={rgb, 255:red, 0; green, 0; blue, 0 }  ][line width=0.75]    (10.93,-3.29) .. controls (6.95,-1.4) and (3.31,-0.3) .. (0,0) .. controls (3.31,0.3) and (6.95,1.4) .. (10.93,3.29)   ;
\draw    (371,161.6) -- (649,161.6) ;
\draw [color={rgb, 255:red, 126; green, 211; blue, 33 }  ,draw opacity=1 ]   (423,115.6) .. controls (463,85.6) and (454,75.6) .. (498,69.6) ;
\draw [color={rgb, 255:red, 245; green, 166; blue, 35 }  ,draw opacity=1 ]   (573,115.6) .. controls (585,138.6) and (583,152.6) .. (587,161.6) ;
\draw  [dash pattern={on 4.5pt off 4.5pt}]  (377,115.6) -- (636,115.6) ;
\draw [color={rgb, 255:red, 245; green, 166; blue, 35 }  ,draw opacity=1 ]   (423,115.6) .. controls (410.14,138.6) and (412.29,152.6) .. (408,161.6) ;
\draw [color={rgb, 255:red, 126; green, 211; blue, 33 }  ,draw opacity=1 ]   (573,115.6) .. controls (533,85.6) and (542,75.6) .. (498,69.6) ;
\draw  [dash pattern={on 4.5pt off 4.5pt}]  (498,57.6) -- (498,174.6) ;

\draw (49,119) node [anchor=north west][inner sep=0.75pt]  [color={rgb, 255:red, 208; green, 2; blue, 27 }  ,opacity=1 ] [align=left] {$\displaystyle \gamma _{-}^{low}$};
\draw (98,65) node [anchor=north west][inner sep=0.75pt]  [color={rgb, 255:red, 126; green, 211; blue, 33 }  ,opacity=1 ] [align=left] {$\displaystyle \gamma _{-}^{up}$};
\draw (248,60) node [anchor=north west][inner sep=0.75pt]  [color={rgb, 255:red, 74; green, 144; blue, 226 }  ,opacity=1 ] [align=left] {$\displaystyle \gamma _{+}^{up}$};
\draw (285,123) node [anchor=north west][inner sep=0.75pt]  [color={rgb, 255:red, 245; green, 166; blue, 35 }  ,opacity=1 ] [align=left] {$\displaystyle \gamma _{+}^{low}$};
\draw (157,120) node [anchor=north west][inner sep=0.75pt]   [align=left] {$\displaystyle y=s_{0}$};
\draw (387,123) node [anchor=north west][inner sep=0.75pt]  [color={rgb, 255:red, 208; green, 2; blue, 27 }  ,opacity=1 ] [align=left] {$\displaystyle \gamma _{-}^{low}$};
\draw (427,65) node [anchor=north west][inner sep=0.75pt]  [color={rgb, 255:red, 126; green, 211; blue, 33 }  ,opacity=1 ] [align=left] {$\displaystyle \gamma _{-}^{up}$};
\draw (560,72) node [anchor=north west][inner sep=0.75pt]  [color={rgb, 255:red, 74; green, 144; blue, 226 }  ,opacity=1 ] [align=left] {$\displaystyle \gamma _{+}^{up}$};
\draw (593,124) node [anchor=north west][inner sep=0.75pt]  [color={rgb, 255:red, 245; green, 166; blue, 35 }  ,opacity=1 ] [align=left] {$\displaystyle \gamma _{+}^{low}$};
\draw (476,39) node [anchor=north west][inner sep=0.75pt]  [font=\small] [align=left] {$\displaystyle HVP( \gamma )$};

\end{tikzpicture}
\caption{Representation of the symmetrization by sections process we have just described for case 2b) ($\gamma$ at the left and $\gamma^r$, with vertical symmetry, at the right).}
\label{symmetry}
\end{figure}

\subsection{Proof of Theorem \ref{theorem A}}
\begin{proof}[Proof of Theorem \ref{theorem A}]
By scaling, we may assume $\lambda=1$. We begin by drawing $HHP(\gamma)$ and denoting by $\gamma_-$ and $\gamma_+$ the portions of the curve lying below and above this line, respectively. Observe that $\gamma_-$ is not necessarily connected. 

By Proposition \ref{symmetry proposition} and the curved polygon approximation, we may assume that both $\gamma_-$ and $\gamma_+$ present vertical symmetry and have piecewise constant curvature. Furthermore, since differentiability is not required for the hinge argument, we may discard the small circle regularization introduced in Section \ref{subsubsection: small circles approxiamtion}. Hence, we can assume that our standing arcs are composed of circular arcs of curvature 1, joined at finitely many non-differentiable points; which we will use to perform the hinge technique (see \cite{diffgeombook} for a proof of the isoperimetric inequality using the hinge technique).

Next, define the angles $\alpha$ and $\beta$ as follows. Let $\alpha$ be the angle between the right tangent to $\gamma$ at the left intersection point of $\gamma$ with $HHP(\gamma)$ and the $HHP(\gamma)$, measured clockwise from the tangent. Similarly, let $\beta$ denote the corresponding angle at the right intersection point, measured counterclockwise from the tangent.

We will first study the case where the contact angles obey $\theta_i\le\pi/2$ for $i=1,2$. Since $\theta_1\le\pi/2$, convexity implies that $\alpha\le\pi/2$; and the same argument applies to $\beta$ since the curve is symmetric. Under these conditions Theorem \ref{theorem B} can be applied and the center of mass of $\gamma_+$ can be minimized by substituting it with the standing circular arc of radius 1 and equal length to $L(\gamma_+)$. 

After this replacement the continuity of the curve may be lost. However, continuity can be restored by vertically translating the segments of $\gamma_-$ until they meet the endpoints of the new circular arc $\gamma_+$.

We name $\alpha'$ and $\beta'$ the angles, defined as $\alpha$ and $\beta$ previously, of the new circle arc obtained after applying Theorem \ref{theorem B} to $\gamma_+$. It is trivial to see that $\alpha'\le\alpha$ and $\beta'\le\beta$.

We now lower the center of mass of $\gamma_-$ using the hinge technique. However, its application needs to be careful, since an unrestricted use of the hinge technique could destroy the continuity of $\gamma$ by turning the components of $\gamma_-$ into individual standing arcs. To avoid this, we preserve a certain constant on the segments $\gamma_-$ during the application of the hinge. 

Let $\{x_i\}_{i=0}^{i=n}$ be the set of non-differentiable points of $\gamma$ between $\gamma(0)$ and $\gamma(L/4)$, ordered from $\gamma(L/4)$ to $\gamma(0)$; with $x_0$ and $x_n$ being $\gamma(L/4)$ and $\gamma(0)$, respectively. For each $i\in[0,n-1]$, let $\vec{u_i}$ be the vector joining $x_{i+1}$ to $x_i$ (which is always contained inside the domain bounded by $\gamma$ due to its convexity). Additionally, from the set of vectors $\{\vec{u_i}\}_{i=0}^{i=n-1}$, we define the set of angles $\{\alpha_i\}_{i=0}^{i=n-1}$, where $\alpha_i$ is the angle between $\vec{u_i}$ and the $x$-axis, measured clockwise from $\vec{u_i}$.

The hinge is now applied to decrease the slopes of each $\vec{u_i}$ ($i\ne 0$), while keeping the direction of $\vec{u_0}$ fixed. The process continues until all non-differentiable points disappear, producing a circular arc of radius 1, equal length and center of mass no greater than the original segment of $\gamma_-$. During this deformation, a downward translation of $\gamma_+$ is performed whenever necessary to preserve continuity. In the following figure, we illustrate an example of the procedure described for an arbitrary segment $\gamma_-$ and $n=3$:

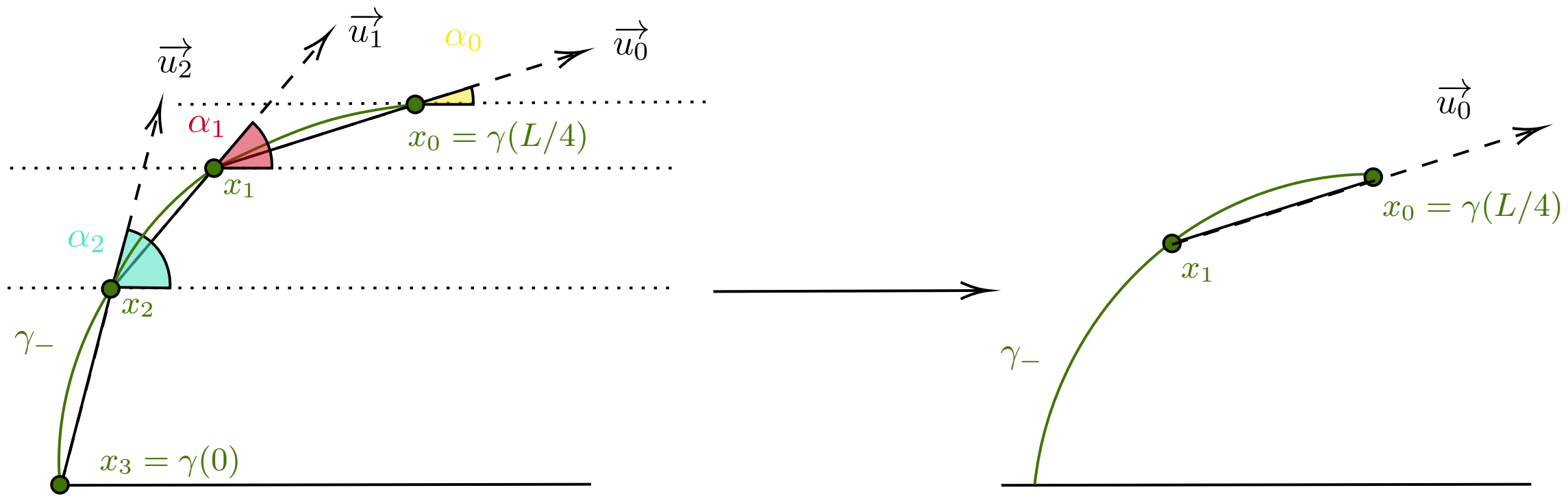
\begin{figure}[h]
    \centering

\tikzset{every picture/.style={line width=0.75pt}} 

\begin{tikzpicture}[x=0.7pt,y=0.7pt,yscale=-1,xscale=1]

\draw    (30,220) -- (238,219.6) ;
\draw    (50,142.6) -- (30,220) ;
\draw    (50,142.6) -- (90,95.6) ;
\draw [color={rgb, 255:red, 65; green, 117; blue, 5 }  ,draw opacity=1 ]   (90,95.6) .. controls (124,79.6) and (135,74.6) .. (169,70.6) ;
\draw    (90,95.6) -- (169,70.6) ;
\draw  [dash pattern={on 4.5pt off 4.5pt}]  (30,220) -- (43.34,168.89) -- (68.49,72.54) ;
\draw [shift={(69,70.6)}, rotate = 104.63] [color={rgb, 255:red, 0; green, 0; blue, 0 }  ][line width=0.75]    (10.93,-3.29) .. controls (6.95,-1.4) and (3.31,-0.3) .. (0,0) .. controls (3.31,0.3) and (6.95,1.4) .. (10.93,3.29)   ;
\draw  [dash pattern={on 4.5pt off 4.5pt}]  (50,142.6) -- (133.7,44.32) ;
\draw [shift={(135,42.8)}, rotate = 130.42] [color={rgb, 255:red, 0; green, 0; blue, 0 }  ][line width=0.75]    (10.93,-3.29) .. controls (6.95,-1.4) and (3.31,-0.3) .. (0,0) .. controls (3.31,0.3) and (6.95,1.4) .. (10.93,3.29)   ;
\draw  [dash pattern={on 4.5pt off 4.5pt}]  (90,95.6) -- (134.03,81.69) -- (233.09,50.4) ;
\draw [shift={(235,49.8)}, rotate = 162.47] [color={rgb, 255:red, 0; green, 0; blue, 0 }  ][line width=0.75]    (10.93,-3.29) .. controls (6.95,-1.4) and (3.31,-0.3) .. (0,0) .. controls (3.31,0.3) and (6.95,1.4) .. (10.93,3.29)   ;
\draw  [dash pattern={on 0.84pt off 2.51pt}]  (9,142.6) -- (269,142.6) ;
\draw  [dash pattern={on 0.84pt off 2.51pt}]  (10,95.6) -- (270,95.6) ;
\draw  [dash pattern={on 0.84pt off 2.51pt}]  (76,70.6) -- (284,69.6) ;
\draw  [fill={rgb, 255:red, 208; green, 2; blue, 27 }  ,fill opacity=0.48 ] (105.2,77.72) .. controls (110.37,82.21) and (113.01,88.83) .. (112.75,95.6) -- (90,95.6) -- cycle ;
\draw  [fill={rgb, 255:red, 248; green, 231; blue, 28 }  ,fill opacity=0.51 ] (191.02,63.82) .. controls (191.59,66.01) and (191.84,68.3) .. (191.75,70.6) -- (169,70.6) -- cycle ;
\draw  [fill={rgb, 255:red, 80; green, 227; blue, 194 }  ,fill opacity=0.58 ] (56.17,119.74) .. controls (58.43,120.3) and (60.62,121.22) .. (62.68,122.54) .. controls (69.59,126.96) and (73.16,134.71) .. (72.75,142.61) -- (50.01,142.32) -- cycle ;
\draw  [draw opacity=0] (29.55,219.78) .. controls (27.83,200.68) and (31.48,179.17) .. (41.09,158.39) .. controls (43.6,152.96) and (46.42,147.8) .. (49.49,142.95) -- (105.03,187.98) -- cycle ; \draw  [color={rgb, 255:red, 65; green, 117; blue, 5 }  ,draw opacity=1 ] (29.55,219.78) .. controls (27.83,200.68) and (31.48,179.17) .. (41.09,158.39) .. controls (43.6,152.96) and (46.42,147.8) .. (49.49,142.95) ;  
\draw  [draw opacity=0] (49.95,142.54) .. controls (54.78,131.24) and (62.87,119.48) .. (73.75,108.94) .. controls (78.47,104.37) and (83.39,100.33) .. (88.4,96.86) -- (100.92,136.99) -- cycle ; \draw  [color={rgb, 255:red, 65; green, 117; blue, 5 }  ,draw opacity=1 ] (49.95,142.54) .. controls (54.78,131.24) and (62.87,119.48) .. (73.75,108.94) .. controls (78.47,104.37) and (83.39,100.33) .. (88.4,96.86) ;  
\draw  [fill={rgb, 255:red, 65; green, 117; blue, 5 }  ,fill opacity=1 ] (26.25,219.78) .. controls (26.25,217.96) and (27.73,216.48) .. (29.55,216.48) .. controls (31.37,216.48) and (32.85,217.96) .. (32.85,219.78) .. controls (32.85,221.61) and (31.37,223.08) .. (29.55,223.08) .. controls (27.73,223.08) and (26.25,221.61) .. (26.25,219.78) -- cycle ;
\draw  [fill={rgb, 255:red, 65; green, 117; blue, 5 }  ,fill opacity=1 ] (46.19,142.95) .. controls (46.19,141.13) and (47.67,139.65) .. (49.49,139.65) .. controls (51.31,139.65) and (52.79,141.13) .. (52.79,142.95) .. controls (52.79,144.78) and (51.31,146.25) .. (49.49,146.25) .. controls (47.67,146.25) and (46.19,144.78) .. (46.19,142.95) -- cycle ;
\draw  [fill={rgb, 255:red, 65; green, 117; blue, 5 }  ,fill opacity=1 ] (86.7,95.6) .. controls (86.7,93.78) and (88.18,92.3) .. (90,92.3) .. controls (91.82,92.3) and (93.3,93.78) .. (93.3,95.6) .. controls (93.3,97.42) and (91.82,98.9) .. (90,98.9) .. controls (88.18,98.9) and (86.7,97.42) .. (86.7,95.6) -- cycle ;
\draw  [fill={rgb, 255:red, 65; green, 117; blue, 5 }  ,fill opacity=1 ] (165.7,70.6) .. controls (165.7,68.78) and (167.18,67.3) .. (169,67.3) .. controls (170.82,67.3) and (172.3,68.78) .. (172.3,70.6) .. controls (172.3,72.42) and (170.82,73.9) .. (169,73.9) .. controls (167.18,73.9) and (165.7,72.42) .. (165.7,70.6) -- cycle ;
\draw    (286,144) -- (392,143.61) ;
\draw [shift={(394,143.6)}, rotate = 179.79] [color={rgb, 255:red, 0; green, 0; blue, 0 }  ][line width=0.75]    (10.93,-3.29) .. controls (6.95,-1.4) and (3.31,-0.3) .. (0,0) .. controls (3.31,0.3) and (6.95,1.4) .. (10.93,3.29)   ;
\draw    (399,220) -- (607,219.6) ;
\draw  [draw opacity=0] (412.1,219.68) .. controls (413.86,203.45) and (418.59,187.18) .. (426.53,171.66) .. controls (450.57,124.72) and (497.5,97.57) .. (545.12,97.95) -- (539.25,229.39) -- cycle ; \draw  [color={rgb, 255:red, 65; green, 117; blue, 5 }  ,draw opacity=1 ] (412.1,219.68) .. controls (413.86,203.45) and (418.59,187.18) .. (426.53,171.66) .. controls (450.57,124.72) and (497.5,97.57) .. (545.12,97.95) ;  
\draw    (465.98,125.13) -- (544.98,100.13) ;
\draw  [fill={rgb, 255:red, 65; green, 117; blue, 5 }  ,fill opacity=1 ] (541.68,99.13) .. controls (541.68,97.31) and (543.15,95.83) .. (544.98,95.83) .. controls (546.8,95.83) and (548.28,97.31) .. (548.28,99.13) .. controls (548.28,100.95) and (546.8,102.43) .. (544.98,102.43) .. controls (543.15,102.43) and (541.68,100.95) .. (541.68,99.13) -- cycle ;
\draw  [fill={rgb, 255:red, 65; green, 117; blue, 5 }  ,fill opacity=1 ] (462.68,125.13) .. controls (462.68,123.31) and (464.15,121.83) .. (465.98,121.83) .. controls (467.8,121.83) and (469.28,123.31) .. (469.28,125.13) .. controls (469.28,126.95) and (467.8,128.43) .. (465.98,128.43) .. controls (464.15,128.43) and (462.68,126.95) .. (462.68,125.13) -- cycle ;
\draw  [dash pattern={on 4.5pt off 4.5pt}]  (466,125.6) -- (510.03,111.69) -- (609.09,80.4) ;
\draw [shift={(611,79.8)}, rotate = 162.47] [color={rgb, 255:red, 0; green, 0; blue, 0 }  ][line width=0.75]    (10.93,-3.29) .. controls (6.95,-1.4) and (3.31,-0.3) .. (0,0) .. controls (3.31,0.3) and (6.95,1.4) .. (10.93,3.29)   ;

\draw (245,38) node [anchor=north west][inner sep=0.75pt]  [font=\small] [align=left] {$\displaystyle \overrightarrow{u_{0}}$};
\draw (141,33) node [anchor=north west][inner sep=0.75pt]  [font=\small] [align=left] {$\displaystyle \overrightarrow{u_{1}}$};
\draw (66,45) node [anchor=north west][inner sep=0.75pt]  [font=\small] [align=left] {$\displaystyle \overrightarrow{u_{2}}$};
\draw (178.75,39.6) node [anchor=north west][inner sep=0.75pt]  [font=\small,color={rgb, 255:red, 248; green, 231; blue, 28 }  ,opacity=1 ] [align=left] {$\displaystyle \alpha _{0}$};
\draw (78,73.6) node [anchor=north west][inner sep=0.75pt]  [font=\small,color={rgb, 255:red, 208; green, 2; blue, 27 }  ,opacity=1 ] [align=left] {$\displaystyle \alpha _{1}$};
\draw (30.75,118.61) node [anchor=north west][inner sep=0.75pt]  [font=\small,color={rgb, 255:red, 80; green, 227; blue, 194 }  ,opacity=1 ] [align=left] {$\displaystyle \alpha _{2}$};
\draw (10,158) node [anchor=north west][inner sep=0.75pt]  [font=\small,color={rgb, 255:red, 65; green, 117; blue, 5 }  ,opacity=1 ] [align=left] {$\displaystyle \gamma _{-}$};
\draw (164.5,75.7) node [anchor=north west][inner sep=0.75pt]  [font=\footnotesize,color={rgb, 255:red, 65; green, 117; blue, 5 }  ,opacity=1 ] [align=left] {$\displaystyle x_{0} =\gamma ( L/4)$};
\draw (43.5,202.7) node [anchor=north west][inner sep=0.75pt]  [font=\footnotesize,color={rgb, 255:red, 65; green, 117; blue, 5 }  ,opacity=1 ] [align=left] {$\displaystyle x_{3} =\gamma ( 0)$};
\draw (52.01,145.32) node [anchor=north west][inner sep=0.75pt]  [font=\footnotesize,color={rgb, 255:red, 65; green, 117; blue, 5 }  ,opacity=1 ] [align=left] {$\displaystyle x_{2}$};
\draw (92,98.6) node [anchor=north west][inner sep=0.75pt]  [font=\footnotesize,color={rgb, 255:red, 65; green, 117; blue, 5 }  ,opacity=1 ] [align=left] {$\displaystyle x_{1}$};
\draw (568,61) node [anchor=north west][inner sep=0.75pt]  [font=\small] [align=left] {$\displaystyle \overrightarrow{u_{0}}$};
\draw (546.98,103.13) node [anchor=north west][inner sep=0.75pt]  [font=\footnotesize,color={rgb, 255:red, 65; green, 117; blue, 5 }  ,opacity=1 ] [align=left] {$\displaystyle x_{0} =\gamma ( L/4)$};
\draw (467.98,131.43) node [anchor=north west][inner sep=0.75pt]  [font=\footnotesize,color={rgb, 255:red, 65; green, 117; blue, 5 }  ,opacity=1 ] [align=left] {$\displaystyle x_{1}$};
\draw (397,164) node [anchor=north west][inner sep=0.75pt]  [font=\small,color={rgb, 255:red, 65; green, 117; blue, 5 }  ,opacity=1 ] [align=left] {$\displaystyle \gamma _{-}$};

\end{tikzpicture}
\caption{Exemplification on the use of the hinge technique for an arbitrary segment of $\gamma_-$ that leaves constant a certain direction regarding its non-differentiable points.}
\label{constant slope}
\end{figure}

The advantage of this procedure over an unrestricted use of the hinge technique is that the resulting curve remains 1-convex. Indeed, since $\vec{u_0}$ remains unchanged, the angle $\varphi$ between the left-tangent at $\gamma(L/4)$ and the $HHP(\gamma)$, measured clockwise from the tangent, also remains fixed. Recalling from before, after applying Theorem \ref{theorem B} to $\gamma_+$ we had $\alpha'\le\alpha$, and due to the original curve $\gamma$ being convex we have $\varphi\ge\alpha$; hence $\varphi\ge\alpha'$, meaning our new curve maintains the 1-convexity property. 

As an observation, if the segments of $\gamma_-$ did not present any non-differentiable points (they are already circular arcs of radius 1), the hinge procedure is unnecessary and this step may be skipped entirely.

Because of the symmetry of $\gamma$, the same construction applies to the right component of $\gamma_-$. Consequently, we obtain a new curve, still denoted by $\gamma$, with center of mass is no greater than that of the original curve and with $\gamma_-$ and $\gamma_+$ being circular arcs of curvature $1$.

Next, we can perform some reflections on the curve that also lower its center of mass. To start, we draw the segment between $\gamma(0)$ and $\gamma(L/4)$, and since $\gamma$ is still convex, we know the left-most segment of $\gamma_-$ only intersects it at its endpoints. For each point $p$ of $\gamma_-$ between $\gamma(0)$ and $\gamma(L/4)$, draw the line joining $p$ to the midpoint $C$ of the segment just described. Reflect $p$ across $C$ along this line, obtaining a new point $p'$ (see Figure \ref{different kinds of symmetries}). By doing this for every point of both segments of $\gamma_-$ we get a new shape of $\gamma_-$ with a strictly lower center of mass, but that together with $\gamma_+$ is not 1-convex. To correct this, we reflect the new lower arc across the vertical line passing through $\gamma(L/4)$ (see Figure \ref{different kinds of symmetries}). Since vertical reflections preserve the height of the center of mass, the resulting curve remains lower than the original one while recovering the required 1-convexity together with $\gamma_+$ (no need to worry about its union points $\gamma(L/4)$ and $\gamma(3L/4)$, since they have remained invariant throughout the procedure).

\begin{figure}[h]
    \centering

\tikzset{every picture/.style={line width=0.75pt}} 

\begin{tikzpicture}[x=0.75pt,y=0.75pt,yscale=-1,xscale=1]

\draw    (32,219.6) -- (211,220.6) ;
\draw    (247,220.6) -- (426,221.6) ;
\draw    (457,220.6) -- (481.99,220.74) -- (636,221.6) ;
\draw    (206,146) -- (257,146.58) ;
\draw [shift={(259,146.6)}, rotate = 180.65] [color={rgb, 255:red, 0; green, 0; blue, 0 }  ][line width=0.75]    (10.93,-3.29) .. controls (6.95,-1.4) and (3.31,-0.3) .. (0,0) .. controls (3.31,0.3) and (6.95,1.4) .. (10.93,3.29)   ;
\draw    (417,146) -- (468,146.58) ;
\draw [shift={(470,146.6)}, rotate = 180.65] [color={rgb, 255:red, 0; green, 0; blue, 0 }  ][line width=0.75]    (10.93,-3.29) .. controls (6.95,-1.4) and (3.31,-0.3) .. (0,0) .. controls (3.31,0.3) and (6.95,1.4) .. (10.93,3.29)   ;
\draw  [draw opacity=0] (44.03,219.68) .. controls (58.54,174.35) and (103.04,138.78) .. (159.37,127.65) -- (197,251.8) -- cycle ; \draw   (44.03,219.68) .. controls (58.54,174.35) and (103.04,138.78) .. (159.37,127.65) ;  
\draw  [draw opacity=0] (391.97,128.32) .. controls (377.46,173.65) and (332.96,209.22) .. (276.63,220.35) -- (239,96.2) -- cycle ; \draw   (391.97,128.32) .. controls (377.46,173.65) and (332.96,209.22) .. (276.63,220.35) ;  
\draw [color={rgb, 255:red, 74; green, 74; blue, 74 }  ,draw opacity=1 ] [dash pattern={on 0.84pt off 2.51pt}]  (44.03,219.68) -- (159.37,127.65) ;
\draw  [draw opacity=0] (277.03,220.08) .. controls (291.54,174.75) and (336.04,139.18) .. (392.37,128.05) -- (430,252.2) -- cycle ; \draw  [color={rgb, 255:red, 128; green, 128; blue, 128 }  ,draw opacity=1 ] (277.03,220.08) .. controls (291.54,174.75) and (336.04,139.18) .. (392.37,128.05) ;  
\draw [color={rgb, 255:red, 74; green, 74; blue, 74 }  ,draw opacity=1 ] [dash pattern={on 0.84pt off 2.51pt}]  (277.03,220.08) -- (392.37,128.05) ;
\draw  [fill={rgb, 255:red, 208; green, 2; blue, 27 }  ,fill opacity=1 ] (98.5,173.66) .. controls (98.5,171.89) and (99.94,170.46) .. (101.7,170.46) .. controls (103.47,170.46) and (104.9,171.89) .. (104.9,173.66) .. controls (104.9,175.43) and (103.47,176.86) .. (101.7,176.86) .. controls (99.94,176.86) and (98.5,175.43) .. (98.5,173.66) -- cycle ;
\draw  [fill={rgb, 255:red, 208; green, 2; blue, 27 }  ,fill opacity=1 ] (331.5,174.06) .. controls (331.5,172.29) and (332.94,170.86) .. (334.7,170.86) .. controls (336.47,170.86) and (337.9,172.29) .. (337.9,174.06) .. controls (337.9,175.83) and (336.47,177.26) .. (334.7,177.26) .. controls (332.94,177.26) and (331.5,175.83) .. (331.5,174.06) -- cycle ;
\draw [color={rgb, 255:red, 74; green, 144; blue, 226 }  ,draw opacity=1 ] [dash pattern={on 4.5pt off 4.5pt}]  (319,162) -- (351,187.6) ;
\draw [color={rgb, 255:red, 74; green, 144; blue, 226 }  ,draw opacity=1 ] [dash pattern={on 4.5pt off 4.5pt}]  (286,200.2) -- (383,149) ;
\draw [color={rgb, 255:red, 74; green, 144; blue, 226 }  ,draw opacity=1 ] [dash pattern={on 4.5pt off 4.5pt}]  (343,145.4) -- (327,202.4) ;
\draw  [draw opacity=0] (599.53,220.92) .. controls (543.12,205.54) and (501.03,171.16) .. (490.54,129.49) -- (666,109.1) -- cycle ; \draw   (599.53,220.92) .. controls (543.12,205.54) and (501.03,171.16) .. (490.54,129.49) ;  
\draw  [fill={rgb, 255:red, 74; green, 144; blue, 226 }  ,fill opacity=1 ] (315.8,162) .. controls (315.8,160.23) and (317.23,158.8) .. (319,158.8) .. controls (320.77,158.8) and (322.2,160.23) .. (322.2,162) .. controls (322.2,163.77) and (320.77,165.2) .. (319,165.2) .. controls (317.23,165.2) and (315.8,163.77) .. (315.8,162) -- cycle ;
\draw  [fill={rgb, 255:red, 74; green, 144; blue, 226 }  ,fill opacity=1 ] (347.8,187.6) .. controls (347.8,185.83) and (349.23,184.4) .. (351,184.4) .. controls (352.77,184.4) and (354.2,185.83) .. (354.2,187.6) .. controls (354.2,189.37) and (352.77,190.8) .. (351,190.8) .. controls (349.23,190.8) and (347.8,189.37) .. (347.8,187.6) -- cycle ;

\draw (112,175) node [anchor=north west][inner sep=0.75pt]  [font=\small,color={rgb, 255:red, 208; green, 2; blue, 27 }  ,opacity=1 ] [align=left] {$\displaystyle C$};
\draw (301,147) node [anchor=north west][inner sep=0.75pt]  [font=\small,color={rgb, 255:red, 74; green, 144; blue, 226 }  ,opacity=1 ] [align=left] {$\displaystyle p$};
\draw (353,190.6) node [anchor=north west][inner sep=0.75pt]  [font=\small,color={rgb, 255:red, 74; green, 144; blue, 226 }  ,opacity=1 ] [align=left] {$\displaystyle p'$};

\end{tikzpicture}
\caption{Exemplification of the symmetries performed on a segment of $\gamma_-$ that is a circular arc of radius 1}
\label{different kinds of symmetries}
\end{figure}
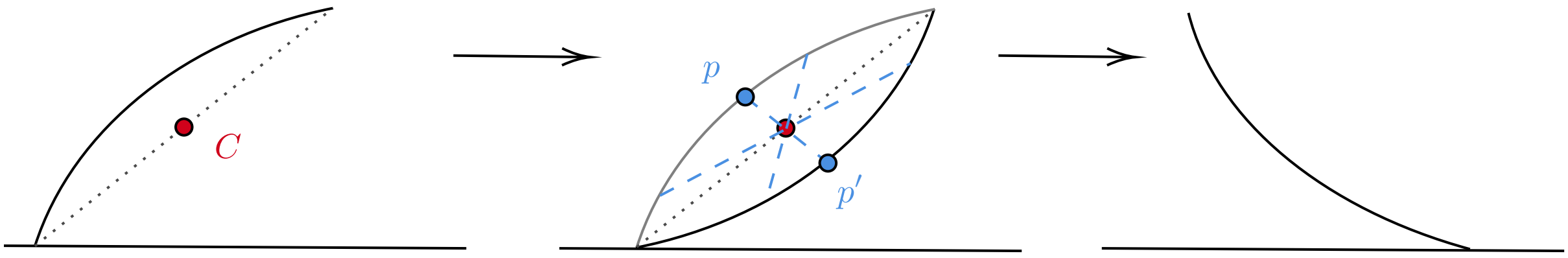

We now regard $HHP(\gamma)$ as a new horizontal reference axis and apply the hinge technique once more, this time attempting to get the center of mass of $\gamma_-$ as close as possible to this line. This will make our curve lose its endpoints with the $x$-axis, allowing us to define $l$ as the lowest horizontal line that intersects $\gamma_-^{\text{new}}$. To get the arc back in contact with the $x$-axis (so that it follows our study's criteria), we will apply a vertical translation until $\gamma_-^{\text{new}}$ has exactly two endpoints with the axis, since the curve we obtain after the hinge technique will keep its symmetry, meaning that both segments of $\gamma_-^{\text{new}}$ will reach the $x$-axis at the same time. We name this new curve $\gamma_-^{\text{low}}$. This can always be done without the segments of $\gamma_-$ intersecting at more than one point, because by the Bent Bow Theorem we have that the standing circular arc has the largest horizontal extension out of all standing convex arcs of piecewise curvature one.

Since we know that these segments are arcs of circle of radius 1 of equal length and with their contact points separated by the exact horizontal distance traveled by the standing arc of circle of radius 1 of length $L(\gamma)/2$ (we know this because $\gamma_+$ has already been transformed into said arc), we can conclude that after this hinge and vertical translation the segments of $\gamma_-$ will intersect exactly at one point: the endpoint with the $x$-axis; obtaining thus the 1-convex lens.

We have to justify that the center of mass has been lowered by this procedure; for that we define the following Euclidean distances:
\begin{gather*}
    d_0=\text{dist}(x-\text{axis}, CM(\gamma_-)), \\
    d_1=\text{dist}(CM(\gamma_-),CM(\gamma_-^{\text{new}})),\\
     d_2=\text{dist}(x-\text{axis},l).
\end{gather*}
Initially,
\begin{equation*}
    CM(\gamma_-)=d_0.
\end{equation*}
After the hinge deformation,
\begin{equation*}
    CM(\gamma_-^{\text{new}})=d_0+d_1.
\end{equation*}
Finally, translating the curve downward until $l$ coincides with the $x$-axis yields
\begin{equation*}
    CM(\gamma^{\text{low}}_-)=d_0+d_1-d_2.
\end{equation*}
By construction, the hinge deformation raises points by at most the distance $d_2$. Using formula \ref{center of mass eqn}, we conclude that $d_1\le d_2$. Therefore,
\begin{equation*}
    CM(\gamma_-^{\text{low}})=d_0+d_1-d_2\le d_0=CM(\gamma_-),
\end{equation*}
which proves that the center of mass has indeed been lowered.

We now consider the case in which the contact angles satisfy $\theta_i>\pi/2$ for $i=1,2$. In this situation $\alpha$ and $\beta$ are not forced to be acute angles, and we must distinguish two separate subcases. 

We first assume that $\alpha,\beta>\pi/2$. Theorem \ref{theorem B} cannot be applied to $\gamma_+$ as done previously. Nevertheless, there exists a section of $\gamma_+$ where Theorem \ref{theorem B} can be applied. More precisely, there exist points $A$,$B\in\gamma$ such that the segment $\overline{AB}$ is parallel to the $x$-axis and such that the angles between the tangents of $\gamma$ at $A$ and $B$ and the line $\overline{AB}$ (measured clockwise from the tangent at $A$ and counterclockwise from the tangent at $B$) are at most $\pi/2$. We choose $A$ and $B$ as low as possible with these properties. In particular, since we have $\alpha,\beta>\pi/2$ the length of the section of $\gamma$ lying above $\overline{AB}$ is strictly greater than 0.

Theorem $B$ can then be applied to the section of $\gamma_+$ lying above the segment $\overline{AB}$; replacing this portion of the curve by the minimizing standing circular arc of radius 1 and equal length. For the portion of $\gamma$ lying below $\overline{AB}$, we apply the hinge technique as described previously for the case $\theta_i\le\pi/2$, deforming the lower arcs until they intersect at the endpoints $\gamma(0)$ and $\gamma(L)$. 

Let $C$ denote this intersection point. After the final downward translation described in the previous case, $C$ becomes the endpoint of our curve. As before, the entire deformation decreases the center of mass of the original curve.

We are left with a vertically symmetric 1-convex isosceles triangle $ABC$, where the vertex $C$ lies on the $x$-axis. Denote the length of the sides as follows: 
\[
\overline{AB}=l_1,
\qquad
\overline{AC}=\overline{BC}=l_2.
\]
We now ask the following optimization question: among all such triangles with fixed perimeter
\begin{equation*}
    l_1+2l_2,
\end{equation*}
which values of $l_1$ and $l_2$ minimize the center of mass?

To answer this, we derive an explicit expression for the center of mass of these triangles and minimize it under the constraints: $l_1+2l_2$ is a fixed value (constant perimeter) and $l_1\le 2l_2$ (which follows from the assumption that the segment $\overline{AB}$ lies above $HHP(\gamma)$).

We first compute the height of the center of mass of the equal sides $AC$ and $CB$ (denoted by $h_2$ in Figure \ref{center of mass final proof}). Using the known expression for the center of mass of a standing circular arc of length at most $\pi$, together with the geometry illustrated in Figure \ref{center of mass final proof}, we obtain formula for the desired center of mass.

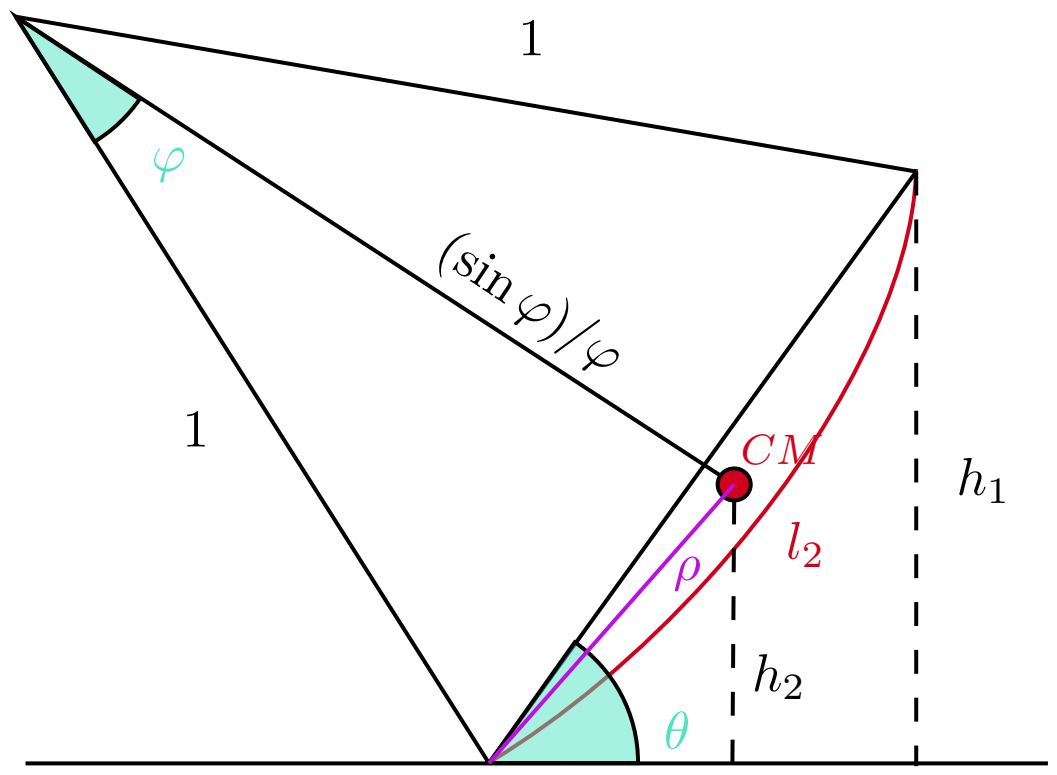
\begin{figure}[h!]
    \centering

\tikzset{every picture/.style={line width=0.75pt}} 

\begin{tikzpicture}[x=0.55pt,y=0.55pt,yscale=-1,xscale=1]

\draw    (116.17,285.17) -- (466,285.17) ;
\draw  [draw opacity=0] (420.87,81.81) .. controls (420.33,119.09) and (394.48,172.64) .. (349.94,221.79) .. controls (326.07,248.12) and (300.2,269.46) .. (275.47,284.43) -- (291.05,171.76) -- cycle ; \draw  [color={rgb, 255:red, 208; green, 2; blue, 27 }  ,draw opacity=1 ] (420.87,81.81) .. controls (420.33,119.09) and (394.48,172.64) .. (349.94,221.79) .. controls (326.07,248.12) and (300.2,269.46) .. (275.47,284.43) ;  
\draw    (112.84,29.72) -- (274.74,285.08) ;
\draw    (112.84,29.72) -- (420.94,82.7) ;
\draw    (420.94,82.7) -- (274.74,285.08) ;
\draw  [fill={rgb, 255:red, 80; green, 227; blue, 194 }  ,fill opacity=0.5 ] (304.12,243.71) .. controls (317.32,252.9) and (325.82,268.11) .. (325.75,285.09) -- (274.74,285.08) -- cycle ;
\draw    (112.84,29.72) -- (358.65,189.75) ;
\draw  [fill={rgb, 255:red, 208; green, 2; blue, 27 }  ,fill opacity=1 ] (352.98,189.75) .. controls (352.98,186.72) and (355.52,184.26) .. (358.65,184.26) .. controls (361.78,184.26) and (364.31,186.72) .. (364.31,189.75) .. controls (364.31,192.78) and (361.78,195.23) .. (358.65,195.23) .. controls (355.52,195.23) and (352.98,192.78) .. (352.98,189.75) -- cycle ;
\draw [color={rgb, 255:red, 189; green, 16; blue, 224 }  ,draw opacity=1 ]   (358.65,189.75) -- (274.74,285.08) ;
\draw  [dash pattern={on 4.5pt off 4.5pt}]  (358.65,195.23) -- (358.04,285.5) ;
\draw  [dash pattern={on 4.5pt off 4.5pt}]  (420.94,82.7) -- (420.93,286.1) ;
\draw  [fill={rgb, 255:red, 80; green, 227; blue, 194 }  ,fill opacity=0.5 ] (155.23,57.91) .. controls (151.27,63.67) and (146.06,68.67) .. (139.79,72.47) -- (112.84,29.72) -- cycle ;

\draw (167.81,163.01) node [anchor=north west][inner sep=0.75pt]  [font=\small] [align=left] {$\displaystyle 1$};
\draw (282.75,29.16) node [anchor=north west][inner sep=0.75pt]  [font=\small] [align=left] {$\displaystyle 1$};
\draw (332.52,265.46) node [anchor=north west][inner sep=0.75pt]  [font=\small,color={rgb, 255:red, 80; green, 227; blue, 194 }  ,opacity=1 ] [align=left] {$\displaystyle \theta $};
\draw (358.83,170.67) node [anchor=north west][inner sep=0.75pt]  [font=\scriptsize,color={rgb, 255:red, 208; green, 2; blue, 27 }  ,opacity=1 ] [align=left] {$\displaystyle CM$};
\draw (262.94,97.71) node [anchor=north west][inner sep=0.75pt]  [font=\small,rotate=-33.75] [align=left] {$\displaystyle (\sin \varphi ) /\varphi $};
\draw (335.77,212.55) node [anchor=north west][inner sep=0.75pt]  [font=\small,color={rgb, 255:red, 189; green, 16; blue, 224 }  ,opacity=1 ] [align=left] {$\displaystyle \rho $};
\draw (362.68,245.75) node [anchor=north west][inner sep=0.75pt]  [font=\small] [align=left] {$\displaystyle h_{2}$};
\draw (432.66,178.61) node [anchor=north west][inner sep=0.75pt]  [font=\small] [align=left] {$\displaystyle h_{1}$};
\draw (157.19,72.68) node [anchor=north west][inner sep=0.75pt]  [font=\small,color={rgb, 255:red, 80; green, 227; blue, 194 }  ,opacity=1 ] [align=left] {$\displaystyle \varphi $};
\draw (374.35,200.5) node [anchor=north west][inner sep=0.75pt]  [font=\small,color={rgb, 255:red, 208; green, 2; blue, 27 }  ,opacity=1 ] [align=left] {$\displaystyle l_{2}$};

\end{tikzpicture}
\caption{Sketch used for the computation of the height of the center of mass of the sides of equal length of the symmetric 1-convex isosceles triangles $ABC$ ($h_2$), where $C$ is the common vertex of the sides with equal length and is in contact with the $x$-axis}
\label{center of mass final proof}
\end{figure}

With the notation of Figure \ref{center of mass final proof} let: $l_1=2\phi,l_2=2\varphi$.

Then, the center of mass of the symmetric 1-convex isosceles triangle $ABC$ is given by

\begin{equation}
\label{minimize CM}
    CM(l_1,l_2)=\frac{2h_2l_2+l_1\left(2\sin\theta\sin\varphi+\frac{\sin\phi}{\phi}-\cos\phi\right)}{L},
\end{equation}

where

\begin{equation*}
    h_2=\rho\cos\left(\theta-\varphi-\arcsin\left(\frac{\sin^2\varphi}{\rho\varphi}\right)\right),
\end{equation*}

\begin{equation*}
    \theta=\arccos\left(\frac{\sin\phi}{2\sin\varphi}\right),
\end{equation*}

and 

\begin{equation*}
    \rho=\sqrt{1+\frac{\sin^2\varphi}{\varphi^2}-\frac{\sin 2\varphi}{\varphi}}.
\end{equation*}

Minimizing \ref{minimize CM} under the above constraints yields the configuration with the smallest center of mass. This minimization can be carried out analytically (see Appendix \ref{appendix A}). The computation shows that the minimum occurs when
\begin{equation*}
l_1=2l_2,
\end{equation*}
which corresponds to the 1-convex lens. This proves Theorem \ref{theorem A} for the case $\theta_i>\pi/2$ and $\alpha,\beta>\pi/2$.

Lastly, consider the case where $\theta_i>\pi/2$ and $\alpha,\beta\le\pi/2$. In this setting, $\gamma_+$ satisfies the assumptions of Theorem \ref{theorem B}. Hence, its center of mass can be lowered by replacing it with the circular arc of radius 1 and length $L/2$. Let $\alpha'$ and $\beta'$ be the new angles obtained after replacing $\gamma_+$ by the minimizing standing circular arc just as in the case $\theta_i\le\pi/2$. As observed previously, $\alpha'\le\alpha$ and $\beta'\le\beta$.

For $\gamma_-$ we will have to apply the hinge technique by parts. Let $A$ and $B$ denote the intersection points of $\gamma$ with $HHP(\gamma)$, ordered from left to right. Define $\delta$ and $\varepsilon$ as the angles between the left-tangent of $\gamma$ at its left most intersection point with $HHP(\gamma)$ and $HHP(\gamma)$ (measured clockwise from the tangent) and for the right most intersection point (measured counterclockwise), respectively. Because $\gamma$ was already a convex curve, we have that $\delta\ge\alpha'$ and $\varepsilon\ge\beta'$.

Now, define $C$ and $D$ as the points of $\gamma_-$ such that the segment $\overline{CD}$ is parallel to the $x$-axis and such that the tangents of $\gamma$ at $C$ and $D$ are perpendicular to this segment, which always exists under this geometric configuration. By symmetry, it suffices to describe the deformation in the left component of $\gamma_-$.

Consider the portion of $\gamma_-$ between $C$ and $A$. Apply the hinge procedure as described earlier for the case $\theta_i\le\pi/2$; that is, keeping fixed the direction of the segment joining $A$ to the nearest non-differentiable point (recall Figure \ref{constant slope}). In particular, this also leaves the angle $\delta$ invariant. As before, a downward translation of $\gamma_+$ is performed whenever necessary to maintain continuity, and the deformation lowers the center of mass while preserving convexity. 

Let $\phi$ denote the angle between the right-tangent to $\gamma$ at $C$ and the $x$-axis, measured clockwise from the tangent. Convexity implies $\phi\ge\delta\ge\alpha'$. 

We then perform the reflection procedure illustrated in Figure \ref{different kinds of symmetries} on the newly obtained arc between $C$ and $A$.

The next step is to apply the hinge to the portion of $\gamma$ between $\gamma(0)$ and $C$ and with respect to the segment $\overline{CD}$, until the left and right tangents of $\gamma$ at $C$ coincide. And, as we have seen for previous cases, an appropriate vertical translation restores contact with the $x$-axis and ensures that the center of mass of the curve has been lowered further. 

We claim that after this deformation the endpoint of the arc remains $\gamma(0)$. Indeed, after the reflections described above, the angle between the left-tangent of $\gamma$ at $A$ and the $x$-axis (measured counter-clockwise from the tangent) equals $\phi$. Since $\phi\ge\alpha'$ and the arc of $\gamma$ between $A$ and $\gamma(0)$ is a circular arc of radius 1 and length $L/4$, its lowest point remains $\gamma(0)$ (otherwise, we could claim that the highest point of $\gamma_+$ is not $\gamma(L/2)$, which is always false after applying Theorem \ref{theorem B}). Hence, after a final translation, $\gamma(0)$, and equivalently for $\gamma(L)$, still becomes the endpoint of the curve.

At this stage, $\gamma_-$ consists of two circular arcs of radius one and length $L/4$. Either these arcs intersect exactly at their endpoints $\gamma(0)$ and $\gamma(L)$, in which case the curve is already the 1-convex lens, or they do not intersect. In the latter case, we can apply a symmetric hinge deformation to both arcs with respect to $HHP(\gamma)$ and then perform one last downward translation, similar to how we finished the case with $\theta_i(\gamma)\le\pi/2$. This lowers the center of mass and produces the 1-convex lens.

We have therefore shown that, in all possible cases, the curve can be transformed into the $1$-convex lens, without increasing the center of mass. This completes the proof of Theorem \ref{theorem A}.

\end{proof}

\section{Conclusions}

Throughout this paper, we have proposed and solved two original reverse isoperimetric problems concerning the height of the center of mass of $\lambda$-convex standing arcs. We proved that the minimizer of the height of the center of mass among all $\lambda$-convex standing curves of fixed length is the $\lambda$-lens. Furthermore, if we impose an upper bound of $\pi/2$ on the contact angles of the standing arcs, the minimizer becomes the standing circular arc of radius $1/\lambda$.

In order to achieve this, we first reviewed relevant results in convex geometry that have been recurring in the study of many other reverse isoperimetric problems, such as Blaschke's Rolling Theorem, the Bent Bow Theorem, and the approximations by curved polygons and small circles.

To prove our main statements, we additionally established two complementary results regarding the convergence of the minimizers of the heights of the centers of mass of curved polygons under the Hausdorff metric to the desired $\lambda$-convex lens. More significantly, we proved that for any arbitrary $\lambda$-convex standing arc, it is always possible to find a standing arc of equal length and vertical symmetry whose center of mass is no higher than that of the original curve. This result draws a parallel to the first formalizations of the isoperimetric inequality, where Steiner had already shown that the maximizers of the inequality must exhibit total symmetry.

For the proofs of the theorems themselves, we made use of both geometrical tools, such as the arc-length parametrization and the Bent Bow Theorem, along with analytical tools, such as the study of the monotonicity of strategically defined functions and integral signs, in order to draw the desired conclusions.

\subsection{Future outlook}

After proving these results for the two dimensional case, the most natural extension is to see if the problem could be generalized and solved for $\lambda$-convex bodies of arbitrary dimensions.

It would also be interesting to follow a similar line of work to other reverse isoperimetric problems and see if our statement can be extended to model spaces of dimension 2 and of constant non-zero curvature, namely the sphere $\mathbb{S}^2$ and the hyperbolic plane $\mathbb{H}^2$ for the cases of constant positive and negative curvature, respectively. If the minimizer is preserved in these spaces we could take the research even further and explore what can be said about model spaces of arbitrary dimension.

\section*{Acknowledgments}
The author would like to express his gratitude to Dr. Kostiantyn Drach (UB), who has supported this project since its inception during a session of Barcelona's Introduction to Mathematical Research Program (BIMR). His advice and comments have played a crucial role in the development and preparation of this manuscript. The author would also like to thank Dr. Laura Prat (UAB), who co-supervised this project and provided invaluable suggestions regarding several of the arguments used throughout the proofs. Finally, the author wishes to thank Dr. Martí Prats Soler (UAB) for his advice during the final stages of the preparation of this manuscript.

The author acknowledges the use of ChatGPT (OpenAI) as an auxiliary tool in deriving the inequality presented in Appendix \ref{appendix A}. The author takes full responsibility for the mathematical content of the manuscript, including the verification of the inequality and all proofs. The underlying ideas, proof strategies, and mathematical framework were developed by the author, with guidance and suggestions from the researchers acknowledged above.

\bibliographystyle{unsrt}
\bibliography{references}

\appendix
\section{Minimization of equation \eqref{minimize CM}}
\label{appendix A}
 Let \(L\in(0,2\pi)\) (the case $L=0$ is not possible and only geometrical configuration possible for the case $L=2\pi$ is the full circle of radius 1) be fixed. To simplify the notation we will write $\varphi=x$ and rewrite all the formulas as functions of only $x$. Therefore, equation \eqref{minimize CM} can be rewritten as:
\begin{equation}
\label{appendix: CM}
CM(x)=4h(x)x+2\phi(x)
\left(
2\sin(\theta(x))\sin x+
\frac{\sin(\phi(x))}{\phi(x)}
-\cos(\phi(x))
\right),
\end{equation}
where
\begin{equation}
\label{appendix: phi}
\phi(x)=\frac L2-2x,
\end{equation}
\begin{equation}
\label{appendix: h}
h(x)=p(x)\cos\left(
\theta(x)-x-\arcsin\left(
\frac{\sin^2 x}{p(x)x}
\right)
\right),
\end{equation}
\begin{equation}
\label{appendix: theta}
\theta(x)=\arccos\left(\frac{\sin(\phi(x))}{2\sin x}\right),
\end{equation}
and
\begin{equation}
\label{appendix: rho}
\rho(x)=\sqrt{
1+\frac{\sin^2 x}{x^2}-\frac{\sin(2x)}{x}
}.
\end{equation}
Since we have that $0<\phi(x)\le2x$ (the case $\phi(x)=0$ is not compatible with the construction of this particular case of the proof of Theorem \ref{theorem A}), from equation\eqref{appendix: phi} we conclude that $x\in[L/8,L/4)$, which in particular allows us to assume that $0<x<\pi/2$.

Our goal will be to prove that the derivative of \eqref{appendix: CM} is always greater than or equal to 0 for $x\in[L/8,L/4)$ and \(L\in(0,2\pi)\). With that, we will have proven that \eqref{appendix: CM} is an increasing function, meaning its minimum is found at the lowest achievable value of $x$ given our constraints, i.e., $x=L/8$, which is our desired result.

If the context allows it, we will omit the explicit dependence on $x$ of the functions previously defined from the notation going forward.

For simplicity, we write
\[
s:=\sin x,\qquad k:=\cos x,
\]
\[
u:=\sin\phi,\qquad v:=\cos\phi.
\]
We will use the fact that all these variables are positive for $0<\phi, x<\pi/2$ and different from 0 without repeatedly mentioning it.

We also define
\[
A:=\sqrt{4\sin^2 x-\sin^2\phi}=\sqrt{4s^2-u^2}.
\]
From equation \eqref{appendix: theta} we trivially get that
\begin{equation}
\label{appendix: cosine theta}
\cos\theta=\frac{\sin\phi}{2\sin x}=\frac{u}{2s},
\end{equation}
which allows us to write 
\begin{equation}
\label{appendix: sine theta}
\sin\theta=
\sqrt{1-\frac{u^2}{4s^2}}
= \sqrt{\frac{4s^2-u^2}{4s^2}}=
\frac{A}{2s}.
\end{equation}
Next we define:
\[
\alpha:=
\arcsin\left(
\frac{\sin^2 x}{\rho(x)x}
\right)=\arcsin\left(
\frac{s^2}{\rho(x)x}
\right),
\]
which directly gives us
\begin{equation}
\label{appendix: rho sine}
\rho\sin\alpha=\frac{s^2}{x}.
\end{equation}
Additionally, from \eqref{appendix: rho} we can see that
\begin{equation}
\label{appendix: complementary 1}
\rho^2-\frac{s^4}{x^2}
=
1+\frac{s^2}{x^2}-\frac{2sk}{x}
-\frac{s^4}{x^2}
= 1-\frac{2sk}{x}+\frac{s^2(1-s^2)}{x^2}=1-\frac{2sk}{x}+\frac{s^2k^2}{x^2}=
\left(1-\frac{sk}{x}\right)^2.
\end{equation}
Therefore, we get
\[
\rho\cos\alpha=\rho\sqrt{1-\sin^2\alpha}=\sqrt{\rho^2-\rho^2\sin^2\alpha},
\]
which from \eqref{appendix: rho sine} and \eqref{appendix: complementary 1} can be simplified to
\begin{equation}
\label{appendix: rho cosine}
    \rho\cos\alpha=\sqrt{\rho^2-\rho^2\sin^2\alpha}=\sqrt{\rho^2-\frac{s^4}{x^2}}=1-\frac{sk}{x}.
\end{equation}
Using the basic trigonometric formulas of the cosine and sine of the sum of two angles plus the relations \eqref{appendix: rho sine} and \eqref{appendix: rho cosine} it is straightforward to see
\[
\rho\cos(x+\alpha)=k-\frac{s}{x}, \qquad \rho\sin(x+\alpha)=s.
\]
Using this and equation \eqref{appendix: cosine theta} and \eqref{appendix: sine theta}, we can rewrite \eqref{appendix: h} as
\[
h=
\left(k-\frac{s}{x}\right)\cos\theta
+
s\sin\theta=\frac{u}{2s}
\left(k-\frac{s}{x}\right)
+
\frac{A}{2}.
\]
Hence,
\[
4xh
=
2xu\cot x-2u+2xA.
\]
From this last expression and \eqref{appendix: sine theta} we can rewrite the expression of the center of mass \eqref{appendix: CM} as follows
\[
CM(x)
=
2xu\cot x-2u+2xA
+
2\phi
\left(
A+\frac{u}{\phi}-v
\right)=
2x\sin\phi\cot x
+
2(x+\phi)A
-
2\phi\cos\phi
.
\]
It is obvious that studying the sign of the derivative of the function $F(x)=CM(x)/2$ is the equivalent to doing it for $CM(x)$, so we will work with
\[
F(x)=x\sin\phi\cot x+(x+\phi)A-\phi\cos\phi.
\]
Now we derive the expression of $F(x)$. For that, we will use the following derivatives of previously introduced variables:
\[
\phi'(x)=-2, \qquad u'=-2v,
\qquad
v'=2u, \qquad A'
=
\frac{4sk+2uv}{A}.
\]
Therefore,

\[
F'(x)
=
u\frac{k}{s}
-2xv\frac{k}{s}
-xu\frac{1}{s^2}
-A
+(x+\phi)\frac{4sk+2uv}{A}
+2v
-2\phi u.
\]
Now define
\[
q:=\frac{u}{2s},
\qquad
r:=\frac{A}{2s}.
\]
Then
\begin{equation}
\label{appendix: r^2+q^2=1}
q^2+r^2=\frac{u^2+A^2}{4s^2}=\frac{u^2+4s^2-u^2}{4s^2}=1.
\end{equation}
Substituting these expressions into \(F'(x)\), we get
\[
F'(x)
=
2(kq+v-sr)
+
2x\left(
\frac{k+qv}{r}
-
\frac{vk+q}{s}
\right)
+
2\phi\left(
\frac{k+qv}{r}
-
2sq
\right).
\]
Define
\[
H:=kq+v-sr,
\]
\[
D:=r(vk+q)-s(k+qv),
\]
and
\[
J:=k+qv-2srq.
\]
Then $F'(x)$ can be written as
\[
F'(x)
=
2H-\frac{2xD}{rs}+\frac{2\phi J}{r}.
\]
Multiplying by \(As^2\), and using \(A=2sr\), gives
\[
As^2F'(x)
=
4s^2(srH-xD+\phi sJ)\implies AF'(x)=4(srH-xD+\phi sJ).
\]
Additionally, observe that
\[
(kr-sq)H=(kr-sq)(kq+v-sr)
=
k^2rq+krv-ksr^2-skq^2-sqv+s^2qr.
\]
Using \eqref{appendix: r^2+q^2=1} and $k^2+s^2=1$ this last expression can be simplified into
\[
(kr-sq)H
=
krv-sqv+rq-sk.
\]
Therefore,
\begin{equation}
\label{appendix: relation between H and D}
(kr-sq)H
=
r(vk+q)-s(k+qv)=D.
\end{equation}
Thus, the expression for the derivative of $F(x)$ can finally be written as:
\begin{equation}
\label{appendix: derivative final expression}
AF'(x)
=
4
\left[
\frac{s^2q}{k}H
+
(\tan x-x)D
+
\phi sJ
\right].
\end{equation}
This expression is particularly nice since it can be proven that all the terms of the sum on the right side of the equality are greater than or equal to 0 under our constraints. Since \(A>0\), $s^2q/k>0$, $\phi s>0$ and $(\tan x-x)>0$ under our conditions, it suffices to prove that
\[
H\geq 0,
\qquad
D\geq 0,
\qquad
J\geq 0,
\]
in order to conclude that $F'(x)\ge0$.

Before diving into the study of the positivity of $H, D$ and $J$ it will be useful to define the following variables:
\[
y:=\frac{\phi}{2}\implies y:=\frac{L}{4}-x.
\]
\[
m:=\tan x,
\qquad
n:=\tan y,
\]
Since $x\in[L/8,L/4)$, we have that $0<y\le x$. In particular this tells us that $x+y=L/4<\pi/2$, which implies that $\tan(x+y)>0$. From here we can give the following useful inequalities for $n$ and $m$
\begin{equation*}
0< n\leq m,
\end{equation*}
\begin{equation*}
\tan(x+y)=\frac{m+n}{1-mn}\implies1-mn>0.
\end{equation*}
Furthermore we define
\[
S:=\sqrt{1+m^2},
\]
which can be trivially related to previously defined constants as
\begin{equation}
\label{appendix: s i k noves}
s=\frac{m}{S},
\qquad
k=\frac{1}{S},
\end{equation}
and using the expressions for the cosine and sine of the sum of angles
\begin{equation}
\label{appendix: v nova}
v=\cos\phi=\cos(2y)=\cos^2y-\sin^2y=\frac{\cos^2y-\sin^y}{\cos^2y+\sin^2y}=\frac{1-n^2}{1+n^2},
\end{equation}
\begin{equation}
\label{appendix: q nova}
q=\frac{\sin\phi}{2\sin x}
=
\frac{1}{2s}\sin(2y)=\frac{S}{2m}2\cos y\sin y=\frac{S}{2m}\frac{2\cos y\sin y}{\cos^2y+\sin^2y}=\frac{S}{m}\frac{n}{1+n^2}=
\frac{nS}{m(1+n^2)}.
\end{equation}
Lastly, we define
\[
N:=m^2(1+n^2)^2-n^2(1+m^2),
\]
and
\[
W:=\sqrt{N}.
\]
Notice that
\[
N=(m^2-n^2)+m^2n^2(1+n^2)\geq 0,
\]
because \(m\geq n\), and therefore $\sqrt{N}$ is well defined. Then
\begin{equation}
\label{appendix: r nova}
r=\sqrt{1-q^2}=\sqrt{\frac{m^2(1+n^2)^2-n^2S^2}{m^2(1+n^2)^2}}=\sqrt{\frac{m^2(1+n^2)^2-n^2(1+m^2)}{m^2(1+n^2)^2}}=\frac{W}{m(1+n^2)}.
\end{equation}
Using all these relations between variables we can finally prove the positivity of $H, D$ and $J$.
\subsubsection*{1) Positivity of \(H\)}
Recall that
\[
H=kq+v-sr.
\]
Using \eqref{appendix: s i k noves}, \eqref{appendix: v nova}, \eqref{appendix: q nova} and \eqref{appendix: r nova} we can rewrite $H$ as
\[
H=
\frac{n}{m(1+n^2)}
+
\frac{1-n^2}{1+n^2}
-
\frac{W}{S(1+n^2)}=
\frac{1-n^2+\frac nm-\frac WS}{1+n^2}.
\]
Defining
\[
B:=1-n^2+\frac nm,
\]
which is always greater than 0 since $n^2\le nm<1$, we get
\begin{equation}
\label{appendix: positivity H}
H=
\frac{BS-W}{S(1+n^2)}=\frac{BS-W}{S(1+n^2)}\frac{BS+W}{BS+W}=\frac{B^2S^2-W^2}
{S(1+n^2)(BS+W)}.
\end{equation}
The numerator of this last expression can be written as
\begin{equation}
\label{appendix: B^2S^2-W^2}
B^2S^2-W^2
=
(1+m^2)
\left(
1-n^2+\frac nm
\right)^2
-
N=
\frac{(m+n)(1-mn)R}{m^2},
\end{equation}
where
\[
R:=3m^2n-mn^2+m+n.
\]
Since
\[
R=m+n+mn(3m-n)\ge0
\]
since \(0< n\leq m\). Therefore, using $R\ge0$ and $1-mn\ge0$ from equation \eqref{appendix: B^2S^2-W^2} we have that $B^2S^2-W^2\ge0$. Therefore, from here we can go back to equation \eqref{appendix: positivity H} to conclude that
\[
H\geq 0.
\]
\subsubsection*{2) Positivity of \(D\)}
Recall from \eqref{appendix: relation between H and D} that $D$ can be related to $H$ as
\[
D=(kr-sq)H.
\]
Since we have already shown \(H\geq 0\), it remains to prove that
\[
kr-sq\geq 0.
\]
From equations \eqref{appendix: s i k noves}, \eqref{appendix: r nova} and \eqref{appendix: q nova} we see that
\[
kr-sq=
\frac{W}{mS(1+n^2)}-\frac{n}{1+n^2}=\frac{W-mnS}{mS(1+n^2)}.
\]
Rationalizing,
\[
kr-sq
=\frac{W-mnS}{mS(1+n^2)}\frac{W+mnS}{W+mnS}=
\frac{W^2-m^2n^2S^2}
{mS(1+n^2)(W+mnS)}.
\]
Expanding, it can be seen that the numerator of the last expression can be written as
\[
W^2-m^2n^2S^2
=
N-m^2n^2(1+m^2)=
(m-n)(m+n)(1-mn)(1+mn).
\]
Therefore,
\[
kr-sq
=
\frac{
(m-n)(m+n)(1-mn)(1+mn)
}{
mS(1+n^2)(W+mnS)
}.
\]
All factors in the numerator are non-negative, and the denominator is positive. Consequently,
\[
kr-sq\geq 0.
\]
and going back to equation \ref{appendix: relation between H and D} and using $H\ge0$ we get that
\[
D\geq 0.
\]
\subsubsection*{3) Positivity of \(J\)}
Recall that
\[
J=k+qv-2srq.
\]
Using \eqref{appendix: s i k noves}, \eqref{appendix: v nova}, \eqref{appendix: q nova} and \eqref{appendix: r nova} we can rewrite $J$ as
\[
J=
\frac1S
+
\frac{nS(1-n^2)}
{m(1+n^2)^2}
-
\frac{2nW}
{m(1+n^2)^2}.
\]
Putting everything over a common denominator,
\[
J=
\frac{
m(1+n^2)^2
+
n(1-n^2)S^2
-
2nSW
}{
mS(1+n^2)^2
}.
\]
Define
\[
P:=
m(1+n^2)^2
+
n(1-n^2)(1+m^2),
\]
which again is always positive since $1-n^2>0$. Therefore, we can write
\[
J=
\frac{P-2nSW}
{mS(1+n^2)^2}=\frac{P-2nSW}
{mS(1+n^2)^2}\frac{P+2nSW}
{P+2nSW}=
\frac{P^2-4n^2S^2W^2}
{
mS(1+n^2)^2(P+2nSW)
}.
\]
Studying the numerator, it can be seen that it can be rewritten as
\[
P^2-4n^2S^2W^2
=
P^2-4n^2(1+m^2)N=
(m+n)(1+n^2)^2(1-mn)R.
\]
Hence,
\[
J=
\frac{
(m+n)(1+n^2)^2(1-mn)R
}{
mS(1+n^2)^2(P+2nSW)
}=
\frac{
(m+n)(1-mn)R
}{
mS(P+2nSW)
}.
\]
As we have previously mentioned for the other cases, all these terms are positive under our conditions. Therefore, we get
\[
J\geq 0.
\]
We have shown that $H, D$ and $J$ are non-negative. Going back to equation \eqref{appendix: derivative final expression} we recall that this proves that $F'(x)\ge0$. Consequently, since $F'(x)=CM'(x)/2$ we have that
\[
CM'(x)\geq 0.
\]
Thus \(CM(x)\) is increasing for $x\in[L/8,L/4)$. Therefore, the minimum of \(CM(x)\) on this interval is attained at the left endpoint:
\[
x_{\text{min}}=\frac L8,
\]
which proves Theorem \ref{theorem A} for the desired case.
\end{document}